\documentclass{article}
\usepackage{pslatex,pifont}
\usepackage{latexsym,amssymb,amscd}
\usepackage{theorem}
\usepackage{pstricks}
\usepackage{fancyhdr}
\usepackage{a4wide}

\def\longto{\longrightarrow}
\def\coker{\mathrm{coker}}
\def\image{\mathrm{image}}

\def\Gt{\widetilde G}
\def\Ghat{\widehat G}
\def\gt{\widetilde g}

\def\Mt{\widetilde M}
\def\xt{\tilde x}
\def\yt{\tilde y}
\def\zt{\tilde z}
\def\d{\mathrm{d\,}}

\def\Ad{\mathop\mathrm{Ad}\nolimits}

\def\Mhat{\widehat{M}}

\def\H{\mathcal{H}}
\def\Hbar{\overline{\mathcal{H}}}

\def\gg{\mathfrak{g}}
\def\tt{\mathfrak{t}}

\def\JJ{\mathbf{J}}
\def\KK{\mathbf{K}}

\def\RR{\mathbb{R}}
\def\TT{\mathbb{T}}
\def\ZZ{\mathbb{Z}}

\def\half{{\textstyle\frac12}}

\newtheorem{theorem}{Theorem}[section]

\newtheorem{proposition}[theorem]{Proposition}
\newtheorem{corollary}[theorem]{Corollary}

{\theorembodyfont{\normalfont}
\newtheorem{example}[theorem]{Example}
\newtheorem{definition}[theorem]{Definition}
\newtheorem{remark}[theorem]{Remark}

}

\newenvironment{proof}%
       {\addvspace\baselineskip\noindent {\sc Proof:}\quad}%
       {\hfill \ding{114} \par\addvspace\baselineskip}  %$\Box$

\makeatletter \@addtoreset{equation}{section}

\makeatother

\begin{document}

\thispagestyle{plain}

\centerline{\large \bf Symplectic Group Actions and Covering Spaces}

\medskip

\centerline{\textbf{James Montaldi \& Juan-Pablo Ortega}}
\medskip
\centerline{\today}
\bigskip
%\maketitle

\begin{abstract}
For symplectic group actions which are not Hamiltonian there are two
ways to define reduction. Firstly using the cylinder-valued momentum
map and secondly lifting the action to any Hamiltonian cover (such
as the universal cover), and then performing symplectic reduction in
the usual way. We show that provided the action is free and proper,
and the Hamiltonian holonomy associated to the action is closed, the
natural projection from the latter to the former is a symplectic
cover. At the same time we give a classification of all
Hamiltonian covers of a given symplectic group action. The main
properties of the lifting of a group action to a cover are studied.\\[12pt]
\end{abstract}

%%%%%%%%%%%%%%%%%%%%%%%%%%%%%%%%%%%%%%%%%%%%%%%%%%%%%%%%%%%%
\section*{Introduction}

There are many instances of symplectic group actions which are not
Hamiltonian---ie, for which there is no momentum map. These can occur
both in applications \cite{MST03} as well as in fundamental studies
of symplectic geometry \cite{lie group valued maps,B02,G05}.  In
such cases it is possible to define a ``cylinder valued momentum
map''~\cite{CDM}, and then perform symplectic reduction with respect
to this map~\cite{OR06,OR06a}. An alternative approach is to lift
to the universal cover, where the action is always Hamiltonian,
and then to perform ordinary symplectic reduction. The principal purpose of this study is to relate the two procedures. In short we show that under suitable hypotheses, the reduced space obtained from the universal cover is a symplectic cover of the one obtained from the cylinder valued momentum map.

In more detail, suppose a connected Lie group $G$ acts on a connected manifold $M$, and let $N$ be a cover of $M$. Then it may not be possible to lift the action of $G$, but there is a natural lift to universal covers giving an action of $\Gt$ on $\Mt$. This can then be used to define an action of $\Gt$ on the given cover $N$. This general construction is well-known, but we were unable to find its principal properties in the literature, and consequently in Section \ref{Lifting group actions to covering spaces} we establish the main results about these lifted actions. For example, since $N$ can be written as a quotient of $\Mt$ by a subgroup of the group of deck transformations, we use this to determine exactly which subgroup of $\Gt$ acts trivially on $N$. We show that if the action on $M$ is free and proper, then so is the appropriate lifted action on $N$. Further details on such lifted actions (including non-free actions) are available as notes \cite{MO-notes}.

In Section 2 we consider the case where $M$ is a symplectic
manifold, and $G$ acts symplectically on $M$. We consider the covers
of $M$ for which the action is Hamiltonian.  The ``largest'' Hamiltonian
cover of $M$ is of course its universal cover $\Mt$; we give an
explicit expression for its momentum map (Proposition \ref{momentum
map for the lifted g tilde action}) and we use it to define a
subgroup of the fundamental group of $M$ whose corresponding set of
subgroups classifies the Hamiltonian covers
(Corollary~\ref{classification covers Hamiltonian}).  There is also
a \emph{``minimal''} such cover, denoted $\Mhat$ and which was first
introduced in \cite{OR04}, where it is called the \emph{universal
covered space} of $M$; we give here a different interpretation of it
as a quotient of the universal cover.

In Section 3, we consider the cylinder valued momentum map of
\cite{CDM} (where it is defined in a different manner, and called
the {\it ``moment r\'eduit''}). In Theorem \ref{thm:reduced
covering} we see that reduction can be carried out in two equivalent
ways. One can either reduce $M$ with respect to the cylinder valued
momentum map or, alternatively, one can lift the action to the
universal cover $\Mt$ (or on any other Hamiltonian cover) and
then carry out (standard) symplectic reduction on  it using its
momentum map. The result is that the natural projection of this
reduced space (inherited from the covering projection) yields the
original reduced space; that is, both reduction schemes are
equivalent up to the projection. If the original action is free and
proper and its Hamiltonian holonomy is closed then both reduced
spaces are symplectic manifolds, and the projection is in fact a
symplectic cover.  We also identify the deck transformation group
of the cover.

We end both sections 2 and 3 with the general example of a group acting by left translations on its cotangent bundle, with symplectic form equal to the sum of the canonical one and a magnetic term consisting of the pullback to the cotangent bundle of a left-invariant 2-form on the group.  In particular we show that symplectic reduction via the cylinder-valued momentum map and Hamiltonian reduction via a standard momentum map yield the same result.

%%%%%%%%%%%%%%%%%%%%%%%%%%
%%%%%%%%%%%%%%%%%%%%%%%%%%
\section{Lifting group actions to covering spaces}
\label{Lifting group actions to covering spaces}

%%%%%%%%%%%%%%%%%%%%%%%%%%
\subsection{The category of covering spaces}\label{sec:covering spaces}
We begin by recalling a few facts about covering spaces. Many of the
details can be found in any introductory book on Algebraic Topology,
for example Hatcher \cite{H}. Let $(M,z_0)$ be a connected manifold
with a chosen base point $z_0$, and let $q_{M}:(\Mt,\zt_0)\to(M,z_0)$ be the universal cover. We realize the universal cover as the set of homotopy classes of paths in $M$ with base point $z_0$. For definiteness, we take the base point in $\Mt$ to be the homotopy class $\zt_0$ of the trivial loop at $z_0$.
Throughout, `homotopic paths' will mean homotopy with fixed
end-points, all paths will be parametrized by $t\in[0,1]$, and for
composition of paths $a*b$ means first do $a$ and then $b$.

Any cover $p_{N}:(N,y_0)\to(M,z_0)$ has the same universal cover
$(\Mt,\zt_0)$ as $(M,z_0)$, and the covering map
$q_{N}:(\Mt,\zt_0)\to(N,y_0)$ can be constructed as follows: Let
$\zt\in \Mt$ and let $z(t)$ be a representative path in $M$, so
$z(0)=z_0$. By the path lifting property of the covering map
$p_{N}$, $z(t)$ can be lifted uniquely to a path $y(t)$ in
$(N,y_0)$. Then $q_{N}(\zt) = y(1)$.

Let $\mathfrak{C}$ be the category of all covers of $(M,z_0)$. The
morphisms are the covering maps. Since any element
$(N,y_0)\in\mathfrak{C}$ also shares $\Mt$ as universal cover, it
sits in a diagram,
$$(\Mt,\zt_0) \stackrel{q_N}{\longrightarrow} (N,y_0)
 \stackrel{p_N}{\longrightarrow} (M,z_0).
$$
Note that the map $\Mt\to M$ can be written both as $q_M$ and as $ p_{\Mt}$.

It is well-known that this category is isomorphic to the category of
subgroups of the fundamental group $\pi_1(M,z_0)$ of $M$, where the
morphisms are the inclusion homomorphisms of subgroups.  The
isomorphism is defined as follows. Let $p_{N}:(N,y_0)\to(M,z_0)$ be
a cover. Then $\Gamma_N:=p_{N*}(\pi_1(N,y_0))$ is the required
subgroup of $\Gamma:=\pi_1(M,z_0)$. $\Gamma_N$ consists of the
homotopy classes of closed paths in $(M,z_0)$ whose lift to
$(N,y_0)$ is also closed, and the number of sheets of the cover
$p_N$ is equal to the index $\Gamma:\Gamma_N$. Note that since $\Mt$
is simply connected, $\Gamma_{\Mt}$ is trivial.

The inverse of this isomorphism can be defined using deck
transformations. Let $\Gamma=\pi_1(M,z_0)$. Then $\Gamma$ is the
fibre of $q_M$ over $z_0$, and it acts on $\Mt$ by deck
transformations defined via the homotopy product: if
$\gamma\in\Gamma$ and $\zt\in\Mt$ then $\gamma*\zt$ gives the action
of $\gamma$ on $\zt$. Then given $\Gamma_1<\Gamma$, define
$N=\Mt/\Gamma_1$, and put $y_0=\Gamma_1\zt_0$. Then from the long
exact sequence of homotopy, it follows that
$\pi_1(N,y_0)\simeq\Gamma_1$. Furthermore, if
$\Gamma_1<\Gamma_2<\Gamma$ then there is a well-defined morphism
(covering map) $p:N_1\to N_2$, where $N_j=\Mt/\Gamma_j$, obtained
from noting that any $\Gamma_1$-orbit is contained in a unique
$\Gamma_2$-orbit, so we put $p(\Gamma_1\zt) = \Gamma_2\zt$.

Let $(N_1,y_1)$ be a cover of $(M,z_0)$ with group $\Gamma_1$, and
let $\Gamma_2=\gamma\Gamma_1\gamma^{-1}$ be a subgroup conjugate to
$\Gamma_1$ (where $\gamma\in \Gamma$).  Then $N_2=\Mt/\Gamma_2$ is
diffeomorphic to $N_1$, but the base point is now
$y_2=\Gamma_2\zt_0$.  A diffeomorphism is simply induced from the
diffeomorphism $\zt\mapsto\gamma\cdot\zt$ of $\Mt$ (which does not
in general map $y_1$ to $y_2$).

If $\Gamma_1\lhd\Gamma$ (normal subgroup), then the cover $(N,y_1)$
is said to be a \emph{normal cover}.  In this case the
$\Gamma$-action (by deck transformations) on $\Mt$ descends to an
action on $N$ (with kernel $\Gamma_1$), and $\Gamma/\Gamma_1$ is the
group of deck transformations of the cover $N\to M$.  For a
general cover, the group of deck transformations is isomorphic to
$N_\Gamma(\Gamma_1)/\Gamma_1$, where $N_\Gamma(\Gamma_1)$ is the
normalizer of $\Gamma_1$ in $\Gamma$.  Only for normal covers does
the group of deck transformations act transitively on the sheets of
the cover. See \cite{H} for examples.

Let us emphasize here that we view $\Gamma=\pi_1(M,z_0)$ both as a
group acting on $\Mt$ by deck transformations, and as a discrete
subset of $\Mt$---the fibre over $z_0$. In particular, for
$\gamma\in\Gamma$, $\gamma*\zt_0 = \gamma$.
In other words, $\zt_0$ is the identity element in $\Gamma$.

%%%%%%%%%%%%%%%%%%%%%%%%%%
\subsection{Lifting the group action} \label{sec:lifting}

Now let $G$ be a connected Lie group acting on the connected
manifold $M$, and let $p_N:(N,y_0)\to (M,z_0)$ be a cover. To
define the lifted action on $N$, we first describe the lift to $\Mt$
and then show it induces an action on $N$, using the cover
$q_N:\Mt\to N$.

The action of $G$ on $M$ does not in general lift to an action of
$G$ on $\Mt$ but of the universal cover $\Gt$, which is also defined
using homotopy classes of paths, with base point the identity
element $e$. The covering map is denoted $q_G:\Gt\to G$. So if $\gt$
is represented by a path $g(t)$ then $q_G(\gt) = g(1)$. The product
structure in $\Gt$ is given by pointwise multiplication of paths: if
$\gt_1$ is represented by a path $g_1(t)$ and $\gt_2$ by $g_2(t)$,
then $\gt_1\gt_2$ is represented by the path $t\mapsto
g_1(t)g_2(t)$.

\begin{definition}
Let $\gt\in\Gt$ be represented by a path $g(t)$ (with $g(0)=e$), and
$\zt\in\Mt$ be represented by  a path $z(t)$ (with $z(0)=z_0$). Then
we define $\gt\cdot\zt$ to be $\yt\in \Mt$, where $\yt$ is the
homotopy class represented by the path $t\mapsto g(t)\cdot z(t)$.
It is readily checked that the homotopy class of this path depends
only on the homotopy classes $\gt$ and $\zt$.
\end{definition}

With this definition for the action of $\Gt$ on $\Mt$, it is clear
that the following diagram commutes:
\begin{equation}\label{eq:CD}
\matrix{ \Gt  \times \Mt & \longto & \Mt \cr
        \downarrow && \downarrow \cr
        G \times M & \longto & M  }
\end{equation}
where the vertical arrows are $q_G\times q_M$ and $q_M$
respectively, and the horizontal arrows are the group actions. In
particular,
\begin{equation}\label{eq:project}
\yt = \gt \cdot \zt \quad\Longrightarrow\quad y=g\cdot z
\end{equation}
where for $\zt \in \Mt $ we denote its projection to $M$ by $z$, and
similarly with elements of $\Gt $.

\begin{remark}
A second approach to defining the action of $\Gt$ on $\Mt$ is as
follows. The action of $G$ gives rise to an `action' of the Lie
algebra $\gg$. That is, to each $\xi\in\gg$ there is associated an infinitesimal generator vector field $\xi_M$ on $M$. Let $N\to M$ be any cover. The
covering map is a local diffeomorphism, so the vector fields $\xi_M$
can be lifted to vector fields $\xi_N$ on $N$. Because this covering
map is a local diffeomorphism, this gives rise to an `action' of
$\gg$ on $N$. Now $\gg$ is the Lie algebra of a unique simply
connected Lie group $\Gt$.  To see that the vector fields on $N$ are
complete, so defining an action of $\Gt$, one needs to compare the
local actions on $M$ and $N$. It is not hard to see that the two
definitions of actions of $\Gt$ are equivalent.
\end{remark}

% \begin{lemma}\label{lemma: 3 paths}
% Let $g(t)$ be a path in $G$ with $g(0)=e$, and $z(t)$ a path in $M$
% with $z(0)=z_0$ and $z(1)=z_1$.  Then the following three homotopy
% classes coincide:
% $$g(t)\cdot z(t),\quad [g(t)\cdot z_0]*[g(1)\cdot z(t)],\quad
%  z(t)*[g(t)\cdot z_1],
% $$
% where $*$ is the homotopy product of paths.
% \end{lemma}
%
%The proof is standard  and left as an exercise for the reader.

% \begin{proof}
% Denote the three curves by $a(t), b(t)$ and $c(t)$ respectively. So
% for example,
%   $$c(t) = \cases{z(2t)& if $t\in[0,\frac12]$\cr g(2t-1)\cdot z_1&
% if $t\in[\frac12,1]$}.$$ A homotopy between $a$ and $b$ can be given
% by
%   $$A(t,s)= \cases{g((1+s)t)\cdot z((1-s^2)t)& if $t\leq
%   \frac1{1+s}$\cr g(1)\cdot z((1+s)t-s)& if $t\geq
%   \frac1{1+s}$}.$$
% Then,  $A(t,0)=a(t)$ and $A(t,1)=b(t)$. It is readily checked that
% $A(t,s)$ is continuous. A similar homotopy can be defined between
% $a$ and $c$.
% \end{proof}

% Recall that $\Gamma:=\pi_1(M,z_0)$ acts on $\Mt$ by deck
% transformations; that is, given $\gamma\in\Gamma$ and $\zt\in\Mt$
% then $\gamma\cdot \zt := \gamma*\zt$. This action is transitive on
% fibres of the covering map $q_M$. Furthermore, the fibre
% $q_M^{-1}(z_0)$ is the $\Gamma$-orbit of the constant loop $\zt_0$
% which we identify with $\Gamma$, see equation
% (\ref{eq:identification of pi_1}).

\begin{proposition}\label{prop:commute}
The action of $\Gt$ on $\Mt$ commutes with the deck transformations.
Furthermore, for each $\gt\in\pi_1(G,e)$ the homotopy class
$g(t)\cdot z_0$ lies in the centre of $\pi_1(M,z_0)$.
\end{proposition}

\begin{proof}
First note that if $g(t)$ is a path in $G$ with $g(0)=e$, and $z(t)$ a path in $M$ with $z(0)=z_0$ and $z(1)=z_1$,  then the following three paths are homotopic:
\begin{equation} \label{eq:3 paths}
g(t)\cdot z(t),\quad [g(t)\cdot z_0]*[g(1)\cdot z(t)],\quad
z(t)*[g(t)\cdot z_1].
\end{equation}

Now let $\gt\in\Gt$, $\delta\in\Gamma$ and $\zt\in \Mt$ with
$q_M(\zt)=y\in M$. We want to show that $\gt\cdot(\delta\cdot\zt) =
\delta\cdot(\gt\cdot\zt)$. By (\ref{eq:3 paths}) applied
with $\gamma=\delta*\zt$, we have
$\gt\cdot(\delta\cdot\zt) = [\delta*\zt]*[\gt\cdot y]$,
while again by (\ref{eq:3 paths}) applied with $\gamma=\zt$ we have $\delta\cdot(\gt\cdot\zt) = \delta*[\zt*(\gt\cdot y)]$.
The result follows from the associativity of the homotopy product.

Finally let $\gt\in\pi_1(G,e)$ and $\delta\in\Gamma$. We want to show
that $[\gt\cdot \zt_0]*\delta = \delta*[\gt\cdot \zt_0]$, where
$\zt_0$ is the constant loop at $x$.  By (\ref{eq:3 paths}),
$\delta*[\gt\cdot \zt_0] = \gt\cdot\delta = [\gt\cdot \zt_0]*
\delta$ (since $g(1)=e$), as required.
\end{proof}

Applying this to the left action of $G$ on itself gives the well-known fact that
$\pi_1(G,e)$ lies in the centre of $\Gt$. Consequently the following
is a central extension:
\begin{equation} \label{eq:ses G}
1 \to \pi_1(G,e) \to \Gt  \stackrel{q_G}{\longto} G \to 1.
\end{equation}

Now we are in a position to define the action of $\Gt$ on an
arbitrary cover $(N,y_0)$ of $(M,z_0)$. As in \S\ref{sec:covering
spaces}, let $\Gamma_N=p_{N*}(\pi_1(N,y_0))<\Gamma$. So, $N\simeq
\Mt/\Gamma_N$. That is, a point in $N$ can be identified with a $\Gamma_N$-orbit of points in $\Mt$.

\begin{definition}
\label{definition of the g tilde action on n}
The $\Gt$-action on $N$ is defined simply by
$$\gt\cdot\Gamma_N\zt := \Gamma_N(\gt\cdot\zt).$$
\end{definition}

This is well-defined as the actions of $\Gt$ and $\Gamma$ commute,
by Proposition \ref{prop:commute}.  It is clear too that the
analogues of (\ref{eq:CD}) and (\ref{eq:project}) hold with $N$ in place of $\Mt$.

\begin{proposition} \label{prop:lifted orbits}
Let $p_N:(N,y_0)\to (M,z_0)$ be a covering map. The $\Gt $-orbits on
$N$ are the connected components of the inverse images under $p_N$
of the orbits on $M$. More precisely, if $y\in p_N^{-1}(z)\subset N$
then $\Gt \cdot y$ is the connected component of $p_N^{-1}(G\cdot
z)$ containing $y$. In particular if the $G$-orbits in $M$ are
closed, so too are the $\Gt$-orbits in $N$.
\end{proposition}

\begin{proof}
Let $Z\subset M$ be any submanifold. Then $Z':=p_N^{-1}(Z)$ is a
submanifold of $N$ and the projection $p_N|_{Z'}:Z'\to Z$ is a
cover, and if $Z$ is closed so too is $Z'$.  Moreover, if $Z$ is
$G$-invariant (hence $\Gt$-invariant), then by the equivariance of
$p_N$ so is $Z'$, and if $Z$ is a single orbit, then $Z'$ is a
discrete union of orbits: discrete because $p_N$ is a cover.
Since $\Gt$ is connected, the orbits are the connected components of
$Z'$.\rule{0pt}{10pt}
\end{proof}

%%%%%%%%%%%%%%%%%%%%%%%%%%%%%%%%%%%%%%%%%%%%%%%%
\subsection{The kernel of the lifted action}
\label{sec:kernel of lifted action}

The natural action of $\Gt$ on $\Mt$ described above need not be
effective, even if the action of $G$ on $M$ is, and the kernel is a
subgroup of $\pi_1(G,e)$ which we describe in this section.

Let $\gt\in \pi_1(G,e)$ be represented by a path $g(t)$, with $g(1)=e$. The path $g(t)$ determines an element $[g(t)\cdot z_0]$ in the centre of $\pi_1(M,z_0)$. Moreover, homotopic loops in $G$ give rise to homotopic loops in $M$, so this induces a well-defined homomorphism
\begin{equation} \label{eq:a_z_0}
a_{z_0}:\pi_1(G,e)\to \pi_1(M,z_0),
\end{equation}
whose image lies in the centre of $\pi_1(M,z_0)$, by Proposition \ref{prop:commute}.

\begin{proposition}
\label{prop:action kernel}
\begin{description}
\item [(i)] The kernel $K<\pi_1(G,e)$ of $a_{z_0}$ is independent of $z_0$ and acts trivially on $\Mt$ and hence on every cover of $M$.
\item [(ii)]  If $(N,y_0)$ is a cover of $(M,z_0)$, with associated subgroup
$\Gamma_N$ of $\pi_1(M,z_0)$, then $K_N:=a_{z_0}^{-1}(\Gamma_N)$ is
independent of the choice of base point $y_0$ in $N$, and acts
trivially on $N$.
\item [(iii)]  If $G$ acts effectively on $M$ then $G_N:= \Gt/K_N$ acts
effectively on $N$.
\end{description}
\end{proposition}

Note that since the domain of $a_{z_0}$ is $\pi_1(G,e)$ which is in
the centre of $\Gt$, it follows that $K_N$ is a normal subgroup of
$\Gt$.  And with the notation of the proposition, $K = K_{\Mt}$
since $\Gamma_{\Mt}$ is trivial. We will write $G':= \Gt/K$
for the group acting on $\Mt$.

In particular, if $a_{z_0}$ is trivial then $K=\pi_1(G,e)$ and the
$G$-action on $M$ lifts to an action of $G$ on $\Mt$. That is,
$a_{z_0}$ is the obstruction to lifting the $G$-action. A particular
case is where the action of $G$ on $M$ has a fixed point. If $z_0$
is such a fixed point then $a_{z_0}=0$.  More generally this is true if any (and hence every) $G$-orbit
in $M$ is contractible in $M$, since in that case too $a_{z_0}$ is
trivial. See also Remark \ref{rmk:Gottlieb}

\begin{proof}
{\bf (i)} Let $z_0,z_1\in M$ and let $\eta$ be any path from $z_0$
to $z_1$ (recall we are assuming $M$ is a connected manifold), and
let $\gt\in\pi_1(G,e)$ with a representative path $g(t)$. For
$T\in[0,1]$ define $g^T(t)=g(Tt)$ (for $t\in[0,1]$), so $g^T\in\Gt$.
Then varying $T$ defines a homotopy from $\eta$ to
$(g^T\cdot\zt_0)*(g(T)(\eta))*((g^T)^{-1}\zt_0')$. In particular,
putting $T=1$ shows that $\eta$ is homotopic to
$a_{z_0}(\gt)*\eta*a_{z_1}(\gt^{-1})$, or equivalently that
$$\eta*a_{z_1}(\gt^{-1})*\bar\eta = a_{z_0}(\gt^{-1}),
$$
where $\bar\eta$ is the reverse of the path $\eta$. This composition
of paths defines the standard isomorphism
$\eta_*:\pi_1(M,z_1)\to\pi_1(M,z_0)$. We have shown therefore that
$a_{z_0}=\eta_*\circ a_{z_1}$ , and so both have the same kernel.
That $K$ acts trivially on $\Mt$ follows from the definition of
$a_{z_0}$: let $\zt\in\Mt$ and $\gt\in K$, then $\gt\cdot\zt =
\gt\cdot(\zt_0*\zt) = a_{z_0}(\gt)*\zt = \zt$ (using
(\ref{eq:3 paths})).

\noindent {\bf (ii)} Let $y_0,y_1\in N$, let $z_j=p_N(y_j)\in M$ and
let $\zeta$ be any path from $y_0$ to $y_1$, with $\eta$ its
projection to $M$. The result follows from the fact that the
following diagram commutes (with ${p_{(N,y_j)}}_*$ written $p_{j_*}$):
\begin{center}
\begin{pspicture}(-1,-1.5)(6,1.5)\psset{unit=1.3}
\rput(0,0){$\pi_1(G,e)$}
\rput(2,-1){$\pi_1(M,z_0)$} \rput(2,1){$\pi_1(M,z_1)$}
\rput(4.5,-1){$\pi_1(N,y_0)$} \rput(4.5,1){$\pi_1(N,y_1)$}
\psline{->}(0.4,0.2)(1.6,0.8)  \rput(0.8,0.6){$a_{z_1}$}
\psline{->}(0.4,-0.2)(1.6,-0.8)   \rput(0.8,-0.6){$a_{z_0}$}
\psline{->}(2,0.7)(2,-0.7) \rput(2.2,0){$\eta_*$}
\psline{->}(4.5,0.7)(4.5,-0.7) \rput(4.7,0){$\zeta_*$}
\psline{<-}(2.6,1)(3.9,1) \rput(3.3,0.8){$p_{0_*}$}
\psline{<-}(2.6,-1)(3.9,-1) \rput(3.3,-0.8){$p_{1_*}$}
\end{pspicture}
\end{center}
Writing $N=\Mt/\Gamma_N$, if $\gt\in a_{z_0}^{-1}(\Gamma_N)$ then
$\gt\in K\Gamma_N$ and, $\gt\,\Gamma_N\zt \subset \Gamma_NK\zt =
\Gamma_N\zt$ so $\gt$ acts trivially (using Proposition
\ref{prop:commute} and part (i)).

\noindent {\bf (iii)} Suppose $\gt\in \Gt$ acts trivially on $N$, so for all $y\in
N$, $\gt\cdot y=y$.  Projecting to $M$, this implies that $g(1)\cdot
z = z$ (for all $z\in M$) so $g(1)\in \cap_{z\in M} G_z=\{e\}$. Thus
$\gt\in\pi_1(G,e)$.

To prove the statement, we first consider the case $N=\Mt$. If
$\gt\not\in K$ then $a_{z_0}(\gt)\neq \zt_0 \in \pi_1(M,z_0)$. Since
$\pi_1(M,z_0)$ acts effectively (by deck transformations) on the
fibre $q_M^{-1}(z_0)\simeq\pi_1(M,z_0)\subset \Mt$ it follows that
$a_{z_0}(\gt)$ acts non-trivially, which is in contradiction with
the assumption that $\gt$ acts trivially.

Now suppose $\gt\in\Gt$ acts trivially on $N$.  We have
$\gt\,\Gamma_N\,\zt_0 =\Gamma_N\,\zt_0$, so that $\gt\in\Gamma_N
K=a_{z_0}^{-1}(\Gamma_N)$ as required.
\end{proof}

\begin{proposition} \label{prop:free and proper}
Let $N$ be any cover of $M$. If the action of $G$ on $M$ is free and proper then so is the action of $G_N$ on $N$.
\end{proposition}

\begin{proof}
First suppose $G$ acts freely on $M$, and let $y=\Gamma_N\zt\in p_N^{-1}(z_0) \subset N$. We need to show that the isotropy group $\Gt_y$ for the $\Gt$ action on $N$ is equal to $K_N$. Now, $\gt\cdot y = \gt\Gamma_N\zt = \gt\Gamma_N\gamma\zt_0$, for some $\gamma\in\Gamma$, as $\Gamma$ acts transitively on the fibre over $z_0$ in $\Mt$. So $\gt\cdot y=y$ if and only if, $\gt\Gamma_N\gamma\zt_0 = \Gamma_N\gamma\zt_0$. However, the action of $\gt$ commutes with that of $\Gamma$ so this reduces to $a_{z_0}(\gt)\in\Gamma_N$ as required for the freeness of the $G_N$-action.

To show the $G_N$-action is proper, we need to show that the action map
$\Phi_N:G_N\times N \to N\times N$ is closed and has compact fibres.
The fibre $\Phi_N^{-1}(x,y) = \{(g, y)\in G_N\times N \mid g\cdot x
= y\}$. If this is non-empty, and $h\cdot x = y$ then
$\Phi_N^{-1}(x,y)\simeq h(G_N)_x$, which is a single element of $G_N$ as the action is free.

To see that the action map is closed, consider a sequence
$(g_i,x_i)$ in $G_N\times N$ for which $(g_i\cdot x_i,x_i)$
converges to $(y,z)$. Then of course $x_i\to z$. We claim that
$g_i\cdot z\to y$. This is because,
$$d(g_i\cdot z, y) \leq d(g_i\cdot z, g_i\cdot x_i) + d(g_i\cdot
x_i, y) = d(z, x_i) + d(g_i\cdot x_i, y),
$$
where $d$ is the $G_N$-invariant metric on $N$ defined above. Both
terms on the right tend to 0 so that $d(g_i\cdot z, y)\to0$ as
required.

Now, by Proposition \ref{prop:lifted orbits} the $G_N$-orbits in $N$
are closed and hence there is an $g\in G_N$ with $y=g\cdot z$.  That
is, $g_i\cdot z \to g\cdot z$. Consequently, $g_i(G_N)_z\to
g(G_N)_z$ in $G_N/(G_N)_z$. By taking a slice to the proper
$(G_N)_z$-action on $G$, this can be rewritten as $g_ih_i\to g$ in
$G_N$, for some sequence $h_i\in (G_N)_z$. Since $(G_N)_z$ is
compact, $(h_i)$ has a convergent subsequence, $h_{i_k}\to h$. Then
$g_{i_k}\to gh^{-1}$. It follows therefore that
$(g_{i_k},x_{i_k})\to (gh^{-1},z)$ and $\Phi_N(gh^{-1}, z) = (y,z)$.
\end{proof}

\begin{remark}\label{rmk:Gottlieb}
D.~Gottlieb \cite{Gott65} considered the images in $\pi_1(M,z_0)$ of ``cyclic homotopies'' of a space, which includes the image of $a_{z_0}$ as a particular case.  He showed in particular that $\image(a_{z_0})$ lies in the subgroup $P(M,z_0)$ of $\pi_1(M,z_0)$ consisting of those loops which act trivially on all homotopy groups $\pi_k(M,z_0)$. Furthermore, he showed that if $M$ is homotopic to a compact polyhedron, and the Euler characteristic $\chi(M)\neq0$, then $\image(a_{z_0})=0$, which implies by what we proved above that every group action on such a space lifts (as an action of $G$) to its universal cover.
\end{remark}

%%%%%%%%%%%%%%%%%%%%%%%%%%
\subsection{Orbit spaces and covers for free actions}
It will be useful for Section \ref{sec:Reduction and Hamiltonian
coverings} to compare the orbit spaces $M/G$ and $\Mt/\Gt$ (or
$\Mt/G'$ where $G'=\Gt/K$) when the $G$-action is free and proper,
and more generally with $N/G_N$ when $N$ is a normal cover of $M$.

Let $N$ be a normal cover of $M$ (see the end of \S\ref{sec:covering
spaces}), with associated group $\Gamma_N$. Then there is an action
of $G_N\times\Gamma$ on $N$ (the action of $\Gamma$ by deck
transformations factors through one of $\Gamma/\Gamma_N$, and
commutes with the $G_N$-action, by Proposition \ref{prop:commute}).

\begin{proposition} \label{prop:orbit space covering}
Let $G$ act freely and properly on $M$. Then the natural map
$q_M':\Mt/G' \to M/G$ is a covering map, with deck transformation
group equal to $\coker(a_{z_0})$ acting transitively on the fibres.

More generally, if $p_N:N\to M$ is a normal cover then $p'_N:
N/G_N \to M/G$ is a normal cover with deck transformation group
$\coker(a_{z_0})/\Gamma_N \ \simeq \ \Gamma/(\image(a_{z_0})\Gamma_N)$.
\end{proposition}

\begin{proof}
Since $G$ acts freely and properly on $M$ then $G_N$ acts freely and
properly on $N$, so both $M/G$ and $N/G_N$ are smooth manifolds.
Moreover, since $N$ is a normal cover of $M$, it follows that
$\Delta_N:=\Gamma/\Gamma_N$ acts freely and transitively on the
fibres of the covering map, and so $M \simeq N/\Delta_N$.

Consider the following commutative diagram:
\begin{equation} \label{eq:Gamma x G}
\begin{CD}
\Mt@ >q_N>> N @>p_N>> M \\
 @V\pi_{\Mt}VV @V\pi_{N}VV @V\pi_{M}VV\\
\Mt/G'@ >q'_N>> N/G_N @>p'_N>> M/G \\
\end{CD}
\end{equation}

Since the covers $q_N$ and $p_N$ are local diffeomorphisms, it
follows that slices to the $\Gt$-actions can be chosen in $\Mt$, $N$
and $M$ in a way compatible with the covers. Consequently, the
lower horizontal maps in the diagram are also covers (the
same is true if the cover $N$ is not normal).

First consider the cover $q_M':\Mt/G'\longto M/G$. Since the
action of $\Gamma$ on $\Mt$ commutes with the action of $G'$, it
descends to an action on $\Mt/G'$. Moreover, since $\Mt/\Gamma
\simeq M$, so
$$(\Mt/G')/\Gamma \simeq \Mt/(G'\times\Gamma) \simeq M/G.$$
(All diffeomorphisms $\simeq$ are natural.) Furthermore, since
$\Gamma$ acts transitively on the fibres of $\Mt\to M$, so it does
on the fibres of $\Mt/G'\to M/G$.

We claim that the isotropy subgroup of the action of $\Gamma$ for
any point in $\Mt/G'$ is $\Gamma\,' = \image(a_{z_0})$. Indeed, for the
action of $G'\times\Gamma$ on $\Mt$ the isotropy subgroup of $\xt$
is
$$
 H = \{(\gt,\gamma)\mid \gt\cdot\gamma\cdot\xt=\xt\}.
$$
Clearly then, $(\gt,\gamma)\in H$ implies in particular
$\gt\in\pi_1(G,e)$, and for such $\gt$, $(\gt,\gamma)\cdot\xt =
a_{z_0}(\gt)*\gamma*\xt$ and so $(\gt,\gamma)\in H$ iff
$a_{z_0}(\gt)=\gamma^{-1}$. Thus $\gamma\in\Gamma$ acts trivially on
$\Mt/G'$ if and only if $\exists\gt\in G'$ such that
$a_{z_0}(\gt^{-1})=\gamma$, as required for the claim. Consequently,
for the cover $q_M'$, the deck transformation group is
$\Gamma/\image(a_{z_0}) = \coker(a_{z_0})$, and this acts
transitively on the fibres.

The same argument as above can be used for the more general normal
cover $p_N:N\to M$, with $G'$ replaced by $G_N$ and $\Gamma$ by
$\Gamma/\Gamma_N$.
\end{proof}

\begin{remark}
If $N$ is a cover of $M$ but not a normal cover, then as pointed out
in the proof $N/G$ is still a cover of $M/G$. Moreover, the fibre still
has cardinality $\coker(a_{z_0})/\Gamma_N$, but the latter is not in
this case a group.
\end{remark}

Notice that as $G$ acts freely and properly on $M$, then $\Mt/G'$ is
a connected and simply connected manifold (simply connected because $G'$ is
connected).  Consequently, $\Mt/G'$ is the (a) universal cover of
$M/G$.

%%%%%%%%%%%%%%%%%%%%%%%%%%%%%%%%%%%%%%%%%%%%%%%%%%%%%%%%%%%%
\section{Hamiltonian covers}

For the remainder of the paper, we assume the manifold $M$ is
endowed with a symplectic form $\omega$ and the Lie group $G$ acts
by symplectomorphisms. Notice that any cover $p _N:N \rightarrow  M
$ of  $M$ is also symplectic with form $\omega _N:= p _N^\ast \omega
$ and that, moreover, the lifted action of $\Gt$ (or $G_N$) on $N$
is also symplectic. It follows that the category of all symplectic
covers of $(M,\omega)$ coincides with the category of all
covers of $M$. Furthermore, the deck transformations on $\Mt$ are
also symplectic.

Symplectic Lie group actions are linked at a very fundamental level
with the existence of {\it momentum maps}. Let $\mathfrak{g}$ be the
Lie algebra of $G$ and $\mathfrak{g}^\ast$ its dual. We recall that
a momentum map $\mathbf{J}:M \rightarrow  \mathfrak{g}^\ast$ for the
symplectic $G$-action on $(M, \omega)$ is defined by the condition
that its components $\mathbf{J}^\xi:= \langle \mathbf{J}, \xi\rangle
$, $ \xi \in \mathfrak{g}$, are Hamiltonian functions for  the
infinitesimal generator vector fields $\xi_M (m):=
\left.\frac{d}{dt}\right|_{t=0}\exp t \xi \cdot  m $. The existence
of a momentum map for the action is by no means guaranteed; however,
it could be that the lifted action to a cover has this feature. For
example, if the cover is simply connected (as is $\Mt$), the action
necessarily has a momentum map associated. This remark leads us to
the following definitions.

\begin{definition}
Let $(M , z_0, \omega )$ be  a connected symplectic manifold
endowed with an action of the connected Lie group $G$. We say that
the smooth cover $p _N: (N,y_0) \rightarrow  (M,z_0) $ of
$(M,z_0)$ is a \emph{Hamiltonian cover} of  $(M,z_0 , \omega)$ if
$N$ is connected and the lifted action of $\Gt$ (or $G_N$) on $(N,
\omega_N)$ has a momentum map $\JJ_N : N \rightarrow
\mathfrak{g}^\ast$ associated.
\end{definition}

Note that we keep the base points in the notation as the choice of momentum map depends on the base point.

If the $G$-action on $M$ is already Hamiltonian, then every cover is
naturally a Hamiltonian cover, so the interesting case is where the
symplectic action on $M$ is not Hamiltonian.

The connectedness hypothesis on $N$ assumed in the previous
definition implies that any two momentum maps of the $G _N$-action
on $N$ differ by a constant element in $\mathfrak{g}^\ast$. We will
assume that $\JJ_N$ is chosen so that $\JJ_N(y_0)=0$. (This choice
should perhaps be denoted $\JJ_{(N,y_0)}$, but we will refrain from
the temptation!)

\begin{definition}
\label{definition of the Hamiltonian covering category} Let  $(M,
z_0, \omega)$ be  a connected symplectic manifold and $G $ a
Lie group acting symplectically thereon. Let $\mathfrak{H} $ be the
category whose objects ${\rm Ob}(\mathfrak{H}) $ are the pairs
$$
 (p_N:(N,y_0,\omega_N) \rightarrow (M,z_0,\omega),\; \JJ_N),
$$
where $p _N$ is a Hamiltonian cover of $(M, z_0, \omega) $ and
$\JJ_N:N\to \gg^*$ is the momentum map for the lifted $\Gt$- (or
$G_N$-) action on $N$ satisfying $\JJ_N(y_0)=0$, and whose morphisms
${\rm Mor}(\mathfrak{H}) $ are the smooth maps $p:(N _1,y_1,
\omega_1) \rightarrow (N _2,y_2, \omega_2)$ that satisfy the
following properties:

\begin{description}
\item [(i)] $p$ is a $\Gt$--equivariant symplectic covering map
\item [(ii)]  the following diagram commutes:
\unitlength=4mm
\begin{center}
\begin{picture}(12,10)
\put(0,5){\makebox(0,0){$(N _1,y_1)$}}
\put(10,5){\makebox(0,0){$(N_2,y_2)$}}
\put(5,0.3){\makebox(0,0){$(M,z_0)$}}
\put(5,9.5){\makebox(0,0){$\mathfrak{g}^\ast$}}
\put(1,4){\vector(1,-1){3}} \put(9,4){\vector(-1,-1){3}}
\put(1,6){\vector(1,1){3}} \put(9,6){\vector(-1,1){3}}
\put(2,5){\vector(1,0){6}} \put(1,2){\makebox(0,0){$p _{N _1}$}}
\put(9,2){\makebox(0,0){$p _{N _2}$}}
\put(1,8){\makebox(0,0){$\mathbf{J} _{N _1}$}}
\put(9,8){\makebox(0,0){$\mathbf{J} _{N _2}$}}
\put(5,6){\makebox(0,0){$p$}}
\end{picture}
\end{center}
\end{description}
We will refer to $\mathfrak{H} $ as the category of
\emph{Hamiltonian covers} of $(M,z_0,\omega)$.
\end{definition}

It should be clear that the ingredients $\omega_N$ and $\JJ_N$ are
both uniquely determined by $p_N:(N,y_0)\to (M,z_0)$ (given the
symplectic form on $M$), so $\mathfrak{H}$ is in fact a (full)
subcategory of the category of all covers of $(M,z_0)$.

The category of the Hamiltonian covers of a symplectic manifold
acted upon symplectically by a Lie algebra was studied
in~\cite{OR04}. We will now use the developments in
Section~\ref{Lifting group actions to covering spaces} to recover
those results in the context of group actions. The study that we
carry out in the following paragraphs sheds light on the {\it
universal covered space} introduced in~\cite{OR04} and additionally
will be of much use in Section~\ref{sec:Reduction and Hamiltonian
coverings} where we will spell out in detail the interplay between
Hamiltonian covers and symplectic reduction.

%%%%%%%%%%%%%%%%%%%%%%%%%%
\subsection{The momentum map on the universal cover}

We now start by giving an expression for  the momentum map
associated to the  $\widetilde{G} $-action on  the universal cover
$\widetilde{M} $ of $M$. As far as this momentum map is concerned,
it does not matter if we consider the $\widetilde{G}$ or the $G' $
action (defined after Proposition \ref{prop:action kernel}) since both have the same Lie algebra and the momentum map depends only on the infinitesimal part of the action. Recall that the \emph{Chu map} $\Psi:M \rightarrow Z
^2(\mathfrak{g})$ is defined by
\begin{equation}\label{eq:chu}
\Psi(z)(\xi,\,\eta) := \omega(z)\left(\xi_M(z),\, \eta_M(z)\right).
\end{equation}
for $\xi,\eta\in\gg$.

\begin{proposition}
\label{momentum map for the lifted g tilde action} Let $(M, \omega)$
be a connected symplectic manifold acted upon symplectically by the
connected Lie group $G$. Then, the $\widetilde{G} $-action on
$(\widetilde{M}, \widetilde{\omega}:= q _M^\ast \omega)$ has a
momentum map associated $\mathbf{J}: \widetilde{M} \rightarrow
\mathfrak{g}^\ast$ that can be expressed as follows: realize
$\widetilde{M} $ as the set of homotopy classes of paths in $M$ with
base point $z _0$. Let $ \widetilde{x} \in \widetilde{M} $  and $x
(t) $ an element in the homotopy class $\widetilde{x} $. Then, for
any $\xi \in \mathfrak{g}$
\begin{equation}
\label{eq:momentum map}
\left<\JJ(\widetilde{x}),\;\xi\right> = \int_{[0,1]}
x^*(\mathbf{i}_{\xi_M}\omega) =
\int_0^1 \omega(x (t))\bigl(\xi_M(x(t)),\,\dot x(t)\bigr)\;\d t.
\end{equation}
If\/ $\widetilde{x}\in\pi_1(M,z_0)$ and
$\widetilde{y}\in\widetilde{M}$ then
$\widetilde{x}*\widetilde{y}\in\Mt$ and
\begin{equation} \label{eq:additive}
\mathbf{J}(\widetilde{x}*\widetilde{y}) = \mathbf{J}(\widetilde{x}) + \mathbf{J}(\widetilde{y}).
\end{equation}
The non-equivariance cocycle $\sigma_{\mathbf{J}}: \widetilde{G}
\rightarrow \mathfrak{g}^\ast$ of\/  $\mathbf{J}$ is given by
\begin{equation}
\label{non-equivariance cocycle} \langle
\sigma_{\mathbf{J}}(\widetilde{g}), \xi\rangle = \int_0^1
\Psi(z_0)(\xi_t,\;\eta_t)\;\d t,
\end{equation}
for any $\xi\in  \mathfrak{g}$, $\widetilde{g} \in \widetilde{G} $,
and $g (t)  $ a curve in the homotopy class of $\widetilde{g}$,
where $\xi_t = \mathrm{Ad}_{g(t)^{-1}}\xi$ and $\eta_t =
\left(T_{e}L_{g(t)}\right)^{-1}\dot g(t)$, and $\Psi$ is the Chu map
defined in (\ref{eq:chu}) above.
\end{proposition}

The non-equivariance cocycle is used to define an affine action of $\Gt$ on $\gg^*$ with respect to which the momentum map is equivariant, namely
\begin{equation}\label{eq:affine action}
\gt\cdot\mu = \Ad^*_{g^{-1}}\mu + \sigma_\JJ(\gt).
\end{equation}
Momentum maps are only defined up to a constant; the one
in~(\ref{eq:momentum map}) is normalized to vanish on the trivial
homotopy class $\widetilde{z}_0$ at $z_0$. The
expression~(\ref{eq:momentum map}) is closely related to the one
in~\cite{MW88} for the momentum map of the action of a group $G$ on
the fundamental groupoid of a symplectic $G$-manifold; see Remark
\ref{rmk:groupoid} below.

\begin{proof}
Let $\alpha :={\bf i}_{\xi_M}\omega$.  Since this 1-form on $M$ is
closed, it follows that $\int x^*\alpha$ depends only on the
homotopy class (indeed homology class) of $x$; that is,
$\JJ(\xt)$ is well-defined by (\ref{eq:momentum map}).

To show that that $\mathbf{J} $ is a momentum map for the
$\widetilde{G}$-action on $\widetilde{M}$, we use the Poincar\'e
Lemma on the closed form $\alpha$. Cover the image of $x(t)$ in $M$
by contractible well-chained open sets (open in $M$),
$U_1,\dots,U_n$, with $x(0)=z_0\in U_1$ and $x(1)\in U_n$. We can
enumerate these sets consecutively along the curve $x(t)$, and let
$z_j=x(t_j)\in U_j\cap U_{j+1}$ lie on the curve and $z_n=x(1)$.

On each $U_j$ we can write $\alpha = \d\phi_j$ for some function
$\phi_j$ (in fact a local momentum for $\xi_M$). Then on $U_i\cap
U_j$, $\mu_{i,j}:= \phi_i-\phi_j$ is constant.
Now, with $I=[0,1]$ and $I_j=[t_j,t_{j+1}]$ we have
\begin{equation} \label{eq:mu i j}
\int_I x^*\alpha = \sum_j \int_{I_j}x^*\d\phi_j =
\sum_j(\phi_j(z_{j+1}) - \phi_j(z_{j})) = \phi_n(z_n)-\phi_1(z_0) -
\sum_{j=1}^{n-1}\mu_{j+1,j}\;.
\end{equation}

The covering map $q_M:\Mt\to M,\;\xt\mapsto x(1)$ identifies the
tangent space $T_{\xt}\Mt$ with $T_{x(1)}M$. Let $ \widetilde{v} \in
T_{\xt}\,\Mt$ arbitrary and $v=T_{\widetilde{x}}\,q
_M(\widetilde{v})$. Thus, differentiating (\ref{eq:mu i j}) at $\xt$
in the direction $\widetilde{v}\in T_{\widetilde{x}}\,\Mt$ gives
$$\d\left(\int x^*\alpha\right)(\widetilde{v}) = \d\phi_n(x(1))(v) =
\alpha(x(1))(v) = \omega(\xi_M, v)= \widetilde{\omega}(\xi_{\widetilde{M}}, \widetilde{v}),$$
as required.
The identity~(\ref{eq:additive}) follows from a straightforward
verification.

We conclude by computing the non-equivariance cocycle
$\sigma_{\mathbf{J}}$. By definition, for any $\widetilde{g} \in
\widetilde{G}$ and $\xi \in \mathfrak{g}$
\[
\sigma_{\mathbf{J}}(\widetilde{g})= \mathbf{J}(\widetilde{g}\cdot
\widetilde{x})- \mbox{\rm Ad} ^\ast _{\widetilde{g}^{-1}}\mathbf{J}
(\widetilde{x}),
\]
for any $\widetilde{x}\in \widetilde{M} $. Take $\widetilde{x}=
\widetilde{z} _0$ and use~(\ref{eq:momentum map}). The formula for
$\sigma_{\mathbf{J}}$ then follows by recalling that
$\mathbf{J}(\widetilde{z} _0)=0 $ and that the $G  $-action on $M$
is symplectic.
\end{proof}

\begin{remark} \label{rmk:isotropic}
If the Chu map vanishes at one point, then $\JJ$ is clearly
coadjoint-equivariant. This happens if there is an isotropic orbit
in $M$ (and hence in $\Mt$).
\end{remark}

\begin{remark}\label{rmk:groupoid}
Let $\Pi(M)$ be the fundamental groupoid of $M$, which has a natural
symplectic structure and Hamiltonian action of $G$ derived from
those on $M$, as described by Mikami and Weinstein, \cite{MW88}. The
relationship between the momentum map $\mathcal{J}:\Pi(M)\to\gg^*$
defined in \cite{MW88} and ours is as follows (we thank Rui Loja
Fernandes for explaining this to us). Given the base point $z_0\in
M$ there is a natural cover $\Mt\times\Mt\to\Pi(M)$ (with fibre
$\pi_1(M,z_0)$). The momentum map $\mathcal{J}$ lifts to one on
$\Mt\times\Mt$, and our momentum map is the restriction of this lift
to the first factor $\Mt\times\{\widetilde{z}_0\}$.

Conversely, given our momentum map $\mathbf{J}:\Mt\to \gg^*$, the
map:
$$\Mt\times\Mt\to\gg^*,\quad(\xt,\yt)\mapsto \mathbf{J}(\xt)-\mathbf{J}(\yt)$$
descends to the quotient by $\pi_1(M,z_0)$ and yields the momentum
map $\mathcal{J}:\Pi(M)\to\gg^*$.
\end{remark}

%%%%%%%%%%%%%%%%%%%%%%%%%%
\subsection{The Hamiltonian holonomy and Hamiltonian covers}
\label{sec:holonomy}

\begin{definition}
\label{Hamiltonian holonomy definition} Let $(M, z_0, \omega)$ be a
connected symplectic manifold with symplectic action of the
connected Lie group $G$. Let $\mathbf{J}:\Mt\to\gg^*$ be the
momentum map defined in Proposition~\ref{momentum map for the lifted
g tilde action}. The \emph{Hamiltonian holonomy} $\mathcal{H}$ of
the $G$-action on $(M, \omega )$ is defined as
$\mathcal{H}=\mathbf{J}(\Gamma)$, and for an arbitrary symplectic
cover $p_N:N\to M$, the holonomy group is
$\mathcal{H}_N:=\JJ(\Gamma_N)$, where $\Gamma = \pi_1(M,z_0)$ and
$\Gamma_N = (p_N)_*\bigl(\pi_1(N,y_0)\bigr)$ (as in \S1).
\end{definition}

\begin{proposition}
The symplectic cover $p_N:(N,y_0)\to (M,z_0)$ is Hamiltonian if and
only if $\H_N=0$.
\end{proposition}

\begin{proof}
If the $\Gt$-action on $N$ is Hamiltonian, then the momentum map is
well-defined. This means that if $\gamma$ is any closed loop in $N$,
then $\JJ(\overline\gamma) = 0$, where
$\overline\gamma\in\pi_1(M,z_0)$ is the image under $(p_N)_*$ of the
homotopy class of $\gamma$. Conversely, if $\H_N=0$ then the map
$\JJ : \Mt\to \gg^*$ descends to a map $\JJ_N : \Mt/\Gamma_N \to
\gg^*$, and as described in \S\ref{Lifting group actions to covering
spaces}, $N\simeq \Mt/\Gamma_N$ as covers of $M$.
\end{proof}

Let us emphasize that if $p_N:(N,y_0)\to (M,z_0)$ is a Hamiltonian
cover, then the momentum map $\JJ_N:N\to\gg^*$ is defined uniquely
by the following diagram.
\begin{equation} \label{eq:KK_N}
\begin{CD}
\Mt@>\JJ>>\mathfrak{g}^\ast\\
@Vq_N VV @VV = V\\
N@>\JJ_N >>\mathfrak{g}^\ast
\end{CD}
\end{equation}

As we pointed out in Section~\ref{Lifting group actions to covering
spaces}, the subgroups of the fundamental group $\Gamma=\pi _1(M, z
_0)$ classify the covers of  $M$.  In a similar vein, the following
result shows that the subgroups of the subgroup $\Gamma_0$ of
$\Gamma$ play the same r\^ole with respect to the Hamiltonian covers
of the symplectic $G$-manifold $(M, \omega )$.

Define,
\begin{equation}   \label{eq:Gamma_0}
\Gamma_0 := \mathbf{J}^{-1}(0)\cap q_M^{-1}(z_0) < \pi_1(M,z_0);
\end{equation}
that is, $\Gamma_0 = \ker (\JJ_{|_\Gamma}:\Gamma\to\gg^*)$. It follows that $\Gamma_0\lhd\Gamma$.

\begin{corollary}
\label{classification covers Hamiltonian}
The symplectic cover $p_N:(N,y_0)\to (M,z_0)$ is Hamiltonian if and
only if\/ $\Gamma_N<\Gamma_0$.  Consequently, $\mathfrak{H}$ is
isomorphic to the category of subgroups of\/ $\Gamma_0$.
\end{corollary}

Recall that the category $\mathfrak{S}(\Gamma)$ of subgroups of a
group $\Gamma$ is the category whose objects are the subgroups, and
whose morphisms are the inclusions of one subgroup into another. We
have therefore shown that $\mathfrak{H} \simeq
\mathfrak{S}(\Gamma_0)$. Explicitly, the isomorphism is given by
\begin{equation}
\begin{array}{rcl}
 \mathfrak{H}&\longrightarrow&\mathfrak{S}(\Gamma_0)\\
 \bigl(p_N:(N,y_0)\to (M,z_0),\,\JJ_N\bigr)& \longmapsto & \Gamma_N = (p_N)_*(\pi_1(N,y_0)).
\end{array}
\end{equation}

%%%%%%%%%%%%%%%%%%%%%%%%%%
\subsection{The universal Hamiltonian covering and covered spaces}

As it was shown in the previous section, the Hamiltonian covers
of a symplectic $G$-manifold $(M, \omega)$ are characterized by the
subgroups of $\Gamma_0$.
The cover associated to the smallest possible subgroup, that is,
the trivial group, is obviously the simply connected universal
cover $\widetilde{M} $ of $M$. It is easy to check that this
object satisfies in the category $\mathfrak{H} $ of Hamiltonian
covers, the same universality property that it satisfies in the
general category of covering spaces, that is,   $(p
_{\widetilde{M}}:\widetilde{M} \rightarrow M, \mathbf{J}) \in {\rm
Ob}(\mathfrak{H}) $ and for any other Hamiltonian cover $(p _N:N
\rightarrow M, \JJ_N)$ of $(M, \omega )$ there exists a morphism
$q_{N }:(\widetilde{M}, \widetilde{\omega}) \rightarrow (N,
\omega_N)$ in ${\rm Mor}(\mathfrak{H}) $. Moreover, any other
element in ${\rm Ob}(\mathfrak{H}) $ that has this universality
property is isomorphic to $(p _{\widetilde{M}}:\widetilde{M}
\rightarrow M, \mathbf{J})$ (we have suppressed the dependence on
base points $z_0,y_0,\zt_0$ in this discussion; if they are included
the morphisms become unique---see Remark \ref{rmk:initial and final
objects} below).

A difference between the general category of covering spaces
and the category of Hamiltonian covers arises when we look at the
cover associated to the biggest possible subgroup of $\Gamma _0$,
that is, $\Gamma_0 $ itself. Unlike the situation found for general
covers, where the biggest possible subgroup that one considers is
the fundamental group $\Gamma $ and it is associated to the trivial
(identity) cover, the cover associated to $\Gamma _0 $ is
non-trivial (unless $M$ is already Hamiltonian) and has an
interesting universality property that is ``dual'' to the one
exhibited by the universal cover.
Define $\Mhat := \Mt/\Gamma_0$; it follows from the corollary above
that this Hamiltonian cover is \emph{minimal}.  It was first
introduced under a different guise in \cite{OR04}, where it is called the \emph{universal covered space} of $(M,\omega)$, and defined using a holonomy bundle associated to a flat $\gg^*$-valued connection. Recall from \S\ref{sec:covering spaces} that
a cover $N\to M$ is said to be normal if $\Gamma_N$ is a normal
subgroup of $\Gamma$. Since $\Gamma_0$ is the kernel of a
homomorphism $\Gamma\to\mathcal{H}$, it follows that $\Mhat$ is a
normal cover of $M$. By Proposition \ref{prop:action kernel}, the
group $\widehat{G} := \Gt/a_{z_0}^{-1}(\Gamma_0)$ acts effectively
on $\Mhat$ (as always, we assume that $G$ acts effectively on $M$).

\begin{proposition}
$\Mhat$ is a Hamiltonian normal cover of $M$ with the universal property
that for any given Hamiltonian cover $p_N:N\to M$ of $M$ there is
a Hamiltonian cover $\hat p_N:N\to \Mhat$.
\end{proposition}

\begin{proof}
Since we have shown that $\mathfrak{H}\simeq\mathfrak{S}(\Gamma_0)$,
this property of $\Mhat$ in $\mathfrak{H}$ follows from the
corresponding property of $\Gamma_0$ in $\mathfrak{S}(\Gamma_0)$;
namely that for every subgroup $\Gamma_1$ of $\Gamma_0$ there is an
inclusion  $\Gamma_1\hookrightarrow\Gamma_0$.
\end{proof}

\begin{remark} \label{rmk:initial and final objects}
$(\widetilde{M}, \widetilde{z}_0 )$ and $(\widehat{M}, \hat{z}_0)$
are initial and final objects in the category of Hamiltonian covers
of $(M, z_0)$  with base points; this of course corresponds to the
fact that $1$ and $\Gamma_0$ are initial and final objects in the
category $\mathfrak{S}(\Gamma_0)$.
\end{remark}

%%%%%%%%%%%%%%%%%%%%%%%%%%
\subsection{The connection in $M\times\gg^*$ and a model for the universal covered space}
\label{sec:connection}

The universal covered space $\Mhat$ was introduced in \cite{OR04}
(though there it is denoted $\Mt$) using a connection in
$M\times\gg^*$ proposed in~\cite{CDM}. Here we briefly review that
definition, and show that it is equivalent to the one given above.

Let $(M, \omega )$ be a connected paracompact symplectic manifold
and let $G$ be a connected Lie group that acts symplectically on
$M$. Consider the Cartesian product $M \times \mathfrak{g}^\ast$  and
let $\pi:M\times\mathfrak{g}^\ast\rightarrow M$ be the projection
onto $M$. Consider $\pi$ as the bundle map of the trivial principal
fiber bundle $(M \times \mathfrak{g}^\ast, M, \pi,
\mathfrak{g}^\ast)$ that has $(\mathfrak{g}^\ast,+) $ as Abelian
structure group. The group $(\mathfrak{g}^\ast,+) $ acts on $M
\times \mathfrak{g}^\ast $ by $\nu \cdot (z, \mu):=(z, \mu- \nu)$.
Let $\alpha \in \Omega^1(M \times \mathfrak{g}^\ast;
\mathfrak{g}^\ast)$ be the connection one-form defined by
\begin{equation} \label{definition of alpha connection}
\langle \alpha(z , \mu) (v_z, \nu), \xi\rangle :=
(\mathbf{i}_{\xi_M} \omega) (z) (v _z) -\langle \nu, \xi \rangle,
\end{equation}
where $(z, \mu)\in M \times \mathfrak{g}^\ast $, $(v _z, \nu) \in T
_z M \times \mathfrak{g}^\ast $,  $\langle\cdot , \cdot \rangle $
denotes the natural pairing between $\mathfrak{g}^\ast $ and
$\mathfrak{g}$, and $\xi_M $ is the infinitesimal generator vector
field associated to $\xi\in \mathfrak{g}$.

The connection $\alpha$ is flat.  For $(z_0,0) \in M \times
\mathfrak{g}^\ast$, let $\Mhat\,':=(M \times \mathfrak{g}^\ast)(z_0,
0) $ be the holonomy bundle through $(z_0, 0) $ and let $
\mathcal{H}(z_0 , 0)$ be the holonomy group of $\alpha$ with
reference point $(z_0, 0) $ (which is an Abelian zero dimensional
Lie subgroup of $\mathfrak{g}^\ast$ by the flatness of $\alpha$); in
other words, $\Mhat\,'$ is the maximal integral leaf of the
horizontal distribution associated to $\alpha$ that contains the
point $(z_0,0) $ and it is hence endowed with a natural initial
submanifold structure with respect to $M \times \mathfrak{g}^\ast$.
See for example Kobayashi and Nomizu \cite{KN63} for standard
definitions and properties of flat connections and holonomy bundles.

The principal bundle $(\Mhat\,',M,\hat p, \mathcal{H}) := (\Mhat\,',
M, \pi|_{(M \times \mathfrak{g}^\ast)(z_0, 0)},\mathcal{H}(z_0 , 0))
$ is a reduction of the principal bundle $(M \times
\mathfrak{g}^\ast, M, \pi, \mathfrak{g}^\ast)$. A straightforward
verification shows that $\mathcal{H}(z_0 , 0) $ coincides with the
Hamiltonian holonomy $ {\mathcal H}$ introduced in
Definition~\ref{Hamiltonian holonomy definition}. In this sense, the
momentum map $\mathbf{J} : \widetilde{M} \rightarrow
\mathfrak{g}^\ast$ establishes a relationship between the deck
transformation groups of the universal cover of $M$ and of the
holonomy bundle $ \widehat{p}: \widehat{M}\,' \rightarrow  M $.
Moreover, the holonomy bundle $\widehat{M}\,' $  can be expressed
using $ \mathbf{J} $ as
\begin{equation}
\label{expression of m hat with j} \widehat{M}\,'=\{(q
_M(\widetilde{x}), \mathbf{J} (\widetilde{x}))\mid \widetilde{x} \in
\widetilde{M}\}.
\end{equation}
This expression allows one to check easily that $(\Mhat\,',M,\hat p,
\mathcal{H}) $ is actually a Hamiltonian cover of $M$ with the
symplectic form $ \widehat{\omega}':= \widehat{p} ^\ast \omega $.
The $G_{\widehat{M}\,'}$-action on $\widehat{M}\,' $ is symplectic
and is  induced by the $\widetilde{G}$-action on $\widehat{M}\,' $
given by
\begin{equation}
\label{action on covered}
\widetilde{g} \cdot (x, \mu)=(g \cdot  x,
\mathbf{J}(\widetilde{g} \cdot \widetilde{x}))=(g \cdot  x,
\sigma_{\mathbf{J}}(\widetilde{g})+ \mbox{\rm Ad} ^\ast _{g
^{-1}}\mathbf{J}(\widetilde{x})),
\end{equation}
where $(x, \mu)\in  \widehat{M}\,'$, $g=p_{\widetilde{G}}
(\widetilde{g})$, and $ \widetilde{x} $ is such that
$p_{\widetilde{M}}(\widetilde{x})= x$, and
$\mathbf{J}(\widetilde{x})= \mu $. The $G_{\widehat{M}\,'}$-action
on $\widehat{M}\,' $ has a momentum map
$\widehat{\mathbf{J}}:\widehat{M}\,' \rightarrow  \mathfrak{g}^\ast
$ given by $\widehat{\mathbf{J}}(x, \mu)= \mu $.

\begin{proposition}
The universal covered space $\widehat{M}=\widetilde{M}/ \Gamma_0$ is
symplectomorphic to $\widehat{M}\,'$.
\end{proposition}

\begin{proof}
The required symplectomorphism is implemented by the map
\[
\begin{array}{cccc}
\Theta: & \widetilde{M}/ \Gamma_0 & \longrightarrow & \widehat{M}\,'\\
   &[\widetilde{x}]&\longmapsto & (x (1), \mathbf{J}(\widetilde{x})).
\end{array}
\]
This map is well defined since by~(\ref{eq:additive}), the smooth
map $\theta:  \widetilde{M} \longrightarrow  \widehat{M}\,' $ given
by $\widetilde{x}\longmapsto  (x (1), \mathbf{J}(\widetilde{x}))$ is
$\Gamma_0 $ invariant and hence it drops to the smooth map $\Theta
$. The map $\theta$ is an immersion since for any $v
_{\widetilde{x}}\in T_{\widetilde{x}} \widetilde{M} $ such that $0=T
_ {\widetilde{x}} \theta \cdot v_{\widetilde{x}}= \left( T_
{\widetilde{x}} p_{\widetilde{M}}\cdot  v_{\widetilde{x}},
T_{\widetilde{x}}\mathbf{J} \cdot v_{\widetilde{x}}\right)$, we have
that $T_ {\widetilde{x}} p_{\widetilde{M}}\cdot  v_{\widetilde{x}}=0
$ and hence $v_{\widetilde{x}}=0 $, necessarily. Given that
$\Gamma_0$ is a discrete group, the projection $\widetilde{M}
\rightarrow \widetilde{M}/ \Gamma_0$ is a local diffeomorphism and
hence $\Theta $ is also an immersion. Additionally, by
~(\ref{expression of m hat with j}), the map $\Theta  $ is also
surjective. We conclude by showing that  $\Theta $ is injective. Let
$\widetilde{x}, \widetilde{y} \in \widetilde{M} $ be such that
$\Theta([\widetilde{x}])=\Theta([\widetilde{y}])$. This implies that
\begin{equation}
\label{two implications for injective}
x (1)=y (1) \quad \mbox{ and
that } \quad  \mathbf{J}(\widetilde{x})=\mathbf{J}(\widetilde{y}).
\end{equation}
The first equality in~(\ref{two implications for injective}) implies
that $\widetilde{x}* \widetilde{\overline{y}} \in \pi _1(M, z _0)$,
where $\widetilde{\overline{y}} $ is the homotopy class associated
to the reverse path $\overline{y} $ of $y$. Moreover, by the second
equality in~(\ref{two implications for injective}), it is easy to
check that $\mathbf{J}(\widetilde{x}* \widetilde{\overline{y}})=0 $,
and hence $\widetilde{x}* \widetilde{\overline{y}} \in \Gamma_0 $.
Since $(\widetilde{x}* \widetilde{\overline{y}})* \widetilde{y}=
\widetilde{x}$ we can conclude that
$[\widetilde{x}]=[\widetilde{y}]$, as required. Consequently,
$\Theta$ being a smooth bijective immersion, it is necessarily a
diffeomorphism. A straightforward verification shows that $\Theta \in \mbox{\rm
Mor}(\mathfrak{H})$, which concludes the proof.
\end{proof}

%%%%%%%%%%%%%%%%%%%%%%%%%%
\subsection{Example} \label{sec:example of T^*G}

We apply the ideas developed in this section to the left action of a Lie group $G$ on its cotangent bundle, but with a modified symplectic form.

Let $G$ be a connected Lie group, and let $\theta:\gg\to\gg^*$ be a symplectic cocycle which is not a coboundary, so it represents a non-zero element of $H^1_s(\gg,\gg^*)$ (the subscript meaning \emph{symplectic} cocycles; that is, $\theta$ is skew-symmetric --- see \cite{Souriau} for details). One can also view $\theta$ as a real-valued 2-cocycle $\Sigma:\gg\times\gg\to\RR$ by putting $\Sigma(\xi,\eta) := \left<\theta(\xi),\,\eta\right>$.  Indeed, $H^2(\gg,\RR)\cong H^1_s(\gg,\gg^*)$.

Let $g(t)$ ($t\in[0,1]$) be a differentiable path in $G$ and define,
\begin{equation}\label{eq:Theta cocycle}
\Theta(g(\cdot)) = \int_0^1\Ad_{g(t)^{-1}}^* \theta \left(g(t)^{-1}\dot g(t)\right)\,\d t.
\end{equation}
It is well-known (and easy to check) that $\Theta$ depends only on the homotopy class of the path $g(t)$ (relative to the end points), so by restricting to $g(0)=e$, $\Theta$ defines a map
$\Theta : \Gt \longto \gg^*$.
Moreover, one can also check that $\Theta$ is a 1-cocycle on $\Gt$, and so defines a well-defined element of $H^1(\Gt,\gg^*)$.

Let $\Gamma_0<\pi_1(G,e)$ be the kernel of the restriction of $\Theta$ to the subgroup $\pi_1(G,e)$ of $\Gt$. Then for any subgroup $\Gamma_1<\Gamma_0$, $\Theta$ descends to a 1-cocycle $\Theta_{1}\in H^1(G_1,\,\gg^*)$, where $G_1=\Gt/\Gamma_1$. In particular, write $\widehat{G}=\Gt/\Gamma_0$. (The notation $\Gamma_0$ is justified in the corollary below.)

%\todo{This must all be well-known! So find a reference!}

Now consider the action of $G$ on $T^*G$ by lifting left multiplication. Given the 2-cocycle $\Sigma$ associated to $\theta$, define a closed differential 2-form  $B_\theta$ on $G$ to be the left-invariant 2-form whose value at $e$ is $\Sigma$.
Write $\pi:T^*G \to G$, and on $M=T^*G$ consider the symplectic form
\begin{equation}\label{eq:cocycle momentum map}
\Omega_\theta = \Omega_{\mathrm{canon}} \; - \; \pi^*B_\theta.
\end{equation}
where $\Omega_{\mathrm{canon}}$ is the canonical cotangent bundle symplectic form.

We claim that the action of $G$ on $M$ is symplectic, and is Hamiltonian if and only if $\Gamma_0=\pi_1(G,e)$. More generally, we claim that whenever $\Gamma_1<\Gamma_0$ the lift of the action to $T^*G_1$ is Hamiltonian.

\begin{proposition}
The action of $\Gt$ on $\Mt=T^*\Gt\cong \Gt\times\gg^*$ with symplectic form given by (\ref{eq:cocycle momentum map}) is Hamiltonian, with momentum map given by
$$\JJ^\theta(\gt, \mu) = \Ad^*_{g^{-1}} \mu \;+  \; \Theta(\gt),$$
where $g=\gt(1)$, and we have identified the Lie algebras of $G$ and $\Gt$. The non-equivariance cocycle of this momentum map is simply $\Theta$.
\end{proposition}

If $\theta=\delta\nu$ for some $\nu\in\gg^*$ (ie $\theta$ represents zero in $H^1(\gg,\gg^*)$), then the action on $T^*G$ is Hamiltonian with momentum map $\JJ(g,\mu) = \Ad_{g^{-1}}^*\mu + \nu$.

\begin{proof}
The action is symplectic because $B_\theta$ is left-invariant.  For the momentum map, the first term of the right-hand side in (\ref{eq:cocycle momentum map}) is the standard expression due to $\Omega_{\mathrm{canon}}$. For the second term, one needs to check that
$$-\iota_{\xi_{\Mt}}\pi^*B_\theta = \left<\d\Theta, \xi\right>.$$
Each side of this is an invariant function, so it suffices to check the equality at the identity element. Now, $\iota_{\xi_{\Mt}}\pi^*B_\theta = \iota_{\xi_G}B_\theta$ and at the identity this is $\iota_\xi\Sigma$. On the other hand $\left<\d\Theta(e)(\eta),\xi\right> = \left<\theta(\eta),\xi\right> = -\Sigma(\xi,\eta)$.

For the non-equivariance cocycle $\sigma\in H^1(G,\gg^*)$,
$$\sigma(h) = \JJ^\theta(h\cdot(e,0)) - \Ad_{h^{-1}}^*\JJ^\theta(e,0) = \JJ^\theta(h,0)-0 =  \Theta(h).$$
\vskip-6mm
\end{proof}

Notice that $\JJ^\theta(e,0)=0$, so this choice of momentum map agrees with the one of Proposition \ref{momentum map for the lifted g tilde action} if we take $z_0=(e,0)$ as base point.

\begin{corollary}
The group $\Gamma_0<\pi_1(G,e)$ defined in (\ref{eq:Gamma_0}) coincides with the group $\Gamma_0$ defined above in terms of $\Theta$.  Consequently, given any subgroup $\Gamma_1<\pi_1(G,e)$, the action of $G_1$ on $T^*G_1$ is Hamiltonian if and only if $\Gamma_1<\Gamma_0$.
\end{corollary}

\begin{proof}
Following the notation of \S\ref{sec:holonomy}, we can take $z_0=(e,0)\in M=T^*G$, and $q_M=q_G\times \mathrm{id}$ on $\Mt=T^*\Gt\simeq \Gt\times\gg^*$. Then $q_M^{-1}(z_0) = \pi_1(G,e)\times\{0\}$ and
$$\Gamma_0:=\left(\JJ^\theta\right)^{-1}(0)\cap\left(\pi_1(G,e)\times\{0\}\right) = \Theta^{-1}(0)\cap\pi_1(G,e),
$$
as required. The rest of the statement follows from Corollary \ref{classification covers Hamiltonian}.
\end{proof}

Notice that with $\Ghat = \Gt/\Gamma_0$, $T^*\Ghat$ is the universal covered space for the given symplectic action of $G$, and it depends on the choice of $\theta$.

\begin{example} \label{eg:torus}
Let $G=\TT=\TT^d=\RR^d/\ZZ^d$ be a $d$-dimensional torus, so $\Gt=\RR^d$ and $\pi_1(G,e)=\ZZ^d$, and $\gg=\RR^d$ can be identified with $\Gt$.  For this case, $H^1_s(\tt,\tt^*)$ is the space of all skew-symmetric linear maps $\tt\to\tt^*$. Let $\theta$ be such a map. Then $\Theta:\Gt\to\tt^*$ can be identified with $\theta$, and the subgroup $\Gamma_0<\ZZ^d$ is  $\Gamma_0 = \ker(\theta)\cap \ZZ^d$.
In particular, if $\theta:\tt\to\tt^*$ is invertible then $\Gamma_0=0$ and the only Hamiltonian cover is the universal cover $\RR^d$. The same occurs if $\ker\theta$ is ``sufficiently irrational''.  If, on the other hand, $\ker\theta$ contains some but not all points of the integer lattice, then $\Ghat$ will be a cylinder; that is a product $\TT^r\times\RR^{d-r}$ for some $r$ with $1\leq r \leq d-1$.
The Hamiltonian holonomy is $\H=\theta(\ZZ^d)\subset\tt^*$, which may or may not be closed in $\tt^*$, depending on the ``irrationality'' of $\ker\theta$.  In all cases, the momentum map on the cover $T^*\RR^d$ is given by $\JJ(u,\mu) = \mu+\Theta(u)$.
\end{example}

\begin{example} \label{eg:Heisenberg}
Consider the group $G$ that is a central extension of $\RR^2$ by $S^1$ with cocycle $\half\omega$. That is, as sets $G=S^1\times\RR^2$, with multiplication
\begin{equation}
(\alpha, u)(\beta, v) = (\alpha+\beta + \half\bar\omega(u,v),\, u+v),
\end{equation}
where $\omega$ is the standard symplectic form on $\RR^2$, and $\half\bar\omega(u,v)=\half\omega(u,v) \bmod{1}\in S^1=\RR/\ZZ$. The universal cover of $G$ is the Heisenberg group $H$, with the same multiplication rule but with $\omega$ in place of $\bar\omega$. We identify $\gg$ with $\RR\times\RR^2$, and correspondingly $\gg^*\simeq \RR^*\times(\RR^2)^*$. One finds that
$$H^1_s(\gg,\,\gg^*) \simeq \left\{\left.\pmatrix{0&\sigma\cr -\sigma^T & 0} \right|\; \sigma\in L(\RR^2,\,\RR^*)\right\}.$$
Now fix any non-zero such $\sigma$ and let $\theta$ be the corresponding element of $H^1(\gg,\gg^*)$. The integral of $\theta$ on $H$ given by (\ref{eq:Theta cocycle}) is,
$$\Theta(\alpha,u) = \pmatrix{\sigma(u)\cr -\alpha\,\sigma - \half\sigma(u)\iota_u\omega}.$$
Note that $\Theta$ does not descend to a function on $G$. The momentum map on $T^*H$ is given by
$$\JJ\left((\alpha,u),\pmatrix{\psi\cr \nu}\right) =
 \Ad^*_{(\alpha,u)^{-1}}\pmatrix{\psi\cr \nu} + \Theta(\alpha,u) =
 \pmatrix{\psi + \sigma(u)\cr \nu - \alpha\sigma -(\psi +\half\sigma(u))\iota_u\omega }.
$$
The Hamiltonian holonomy is therefore
$$\H = \JJ(\ZZ,0) = \pmatrix{0\cr \ZZ\,\sigma},$$
which is closed. The cylinder-valued momentum map on $T^*G$ takes values in $C=\gg^*/\H \simeq \RR\times\RR\times S^1$.

\end{example}

We continue these examples at the end of the next section, where we consider symplectic reduction for such actions.

%%%%%%%%%%%%%%%%%%%%%%%%%%%%%%%%%%%%%%%%%%%%%%%%%%%%%%%%%%%%
\section{Symplectic reduction and Hamiltonian covers}
\label{sec:Reduction and Hamiltonian coverings}

Symplectic reduction is a well studied process that prescribes how
to construct symplectic quotients out of the orbit spaces associated
to the symplectic symmetries of a given symplectic manifold. Even
though it is known how to carry this out for fully general
symplectic actions~\cite{OR06}, the implementation of this procedure
is particularly convenient in the presence of a standard momentum
map, that is, when the Hamiltonian holonomy is trivial (this is the
so called symplectic or Meyer-Marsden-Weinstein reduction~\cite{meyer,mwr}).
Unlike the situation encountered in the general case with a
non-trivial Hamiltonian holonomy, the existence of a standard
momentum map implies the existence of a unique canonical symplectic
reduced space. In the light of this remark the notion of Hamiltonian
cover appears as an interesting and useful object for reduction.
More specifically, one may ask whether, given a symplectic action on
a symplectic manifold with non-trivial holonomy and with respect to
which we want to reduce, we could lift the action to a Hamiltonian
cover, perform reduction there with respect to a standard
momentum map, and then project down the resulting space. How would
this compare with the potentially complicated reduction in the
original manifold? The main result in this section shows that indeed
both processes yield essentially the same result.  Furthermore, we
show that this projection down is a cover.

\subsection{The cylinder valued momentum map}

Recall the definition of the holonomy of a symplectic action of $G$
on $M$ given in Definition \ref{Hamiltonian holonomy definition}:
namely, $\mathcal{H} = \JJ(\Gamma)$, where as always,
$\Gamma=\pi_1(M,z_0)$.  Using this definition, equation
(\ref{eq:additive}) can be expressed by saying that $\mathbf{J}$ is
equivariant with respect to $\Gamma$ acting as deck transformations
on $\Mt$ and as translations by elements of $\mathcal{H}$ on
$\gg^*$. It follows that $\mathbf{J}$  descends to another map with
values in $\gg^*/\mathcal{H}$.
% \[
% \begin{CD}
%  \Mt @> \mathbf{J} >> \mathfrak{g}^\ast\\
%  @V q_M VV @VV V\\
%  M@> >> \gg^\ast/ \mathcal{H}.
% \end{CD}
% \]
However, in general this is a difficult object to use as
$\mathcal{H}$ is not necessarily a \emph{closed} subgroup of
$\gg^*$. To circumvent this, we proceed as follows.

Let $\overline{\mathcal{H}}$ be the closure of $\mathcal{H} $ in
$\mathfrak{g}^\ast$. Since $\overline{\mathcal{H}}$ is a closed
subgroup of $(\mathfrak{g}^\ast, +)$, the quotient $C:=
\mathfrak{g}^\ast/ \overline{\mathcal{H}}$ is a cylinder (that is,
it is isomorphic to the Abelian Lie group $\mathbb{R}^a \times
\mathbb{T}^b$ for some $ a, b \in \mathbb{N}$). Let $\pi_C:
\mathfrak{g}^\ast\rightarrow \mathfrak{g}^\ast/\overline{{\mathcal
H}}$ be the projection. Define $\mathbf{K}: M \rightarrow  C $ to be
the map that makes the following diagram commutative:
\begin{equation}
\label{diagram commutative cylinder valued momentum map}
\begin{CD}
\Mt @> \mathbf{J} >> \mathfrak{g}^\ast\\
@V q_M VV @VV\pi_C V\\
M@>\mathbf{K}>> C \hbox to 0pt{$\;=\gg^\ast/ \Hbar$\hss}
\end{CD}
\end{equation}
In other words, $\mathbf{K}$ is defined by  $ \mathbf{K}(z)= \pi_C
(\JJ(\zt))$, where $\zt \in \Mt $ is any path with endpoint $z$. We
will refer to $ \mathbf{K}:M \rightarrow \mathfrak{g}^\ast/
\overline{{\mathcal H}}$ as a \emph{cylinder valued momentum map}
associated to the symplectic $G$-action on $(M, \omega)$. This
object was introduced in~\cite{CDM} using the connection described in \S\ref{sec:connection}, where it is called the {\it ``moment r\'eduit''}.

Any other choice of Hamiltonian cover in place of $\Mt$ would render
the same Hamiltonian holonomy group $\mathcal{H} $ and the same
cylinder valued momentum map. If one chose a different base point
$z_1\in M$ in place of $z_0$ the holonomy group would remain the
same, but the cylinder valued momentum map would differ from $\mathbf{K}$ by a constant in $\mathfrak{g}^\ast/
\Hbar$.

\paragraph{Elementary properties.}
The cylinder valued momentum map is a strict generalization of the
standard (Kostant-Souriau) momentum map since the $G$-action has a
standard momentum map if and only if the holonomy group
$\mathcal{H}$ is trivial. In such a case the cylinder valued
momentum map is a standard momentum map. The cylinder valued
momentum map satisfies Noether's Theorem; that is, for any
$G$-invariant function $h \in C^\infty(M)^{G}$, the flow $F _t $ of
its associated Hamiltonian vector field $X _h$ satisfies the
identity $ \mathbf{K} \circ F _t= \mathbf{K}| _{{\rm Dom}(F _t)}$.
Additionally, using the diagram (\ref{diagram commutative cylinder
valued momentum map}) and identifying $T_zM$ and $T_{\zt}\Mt$ via
$T_{\zt}q_M$, one has that for any $v _z \in T _zM $, $T _z
\mathbf{K} ( v _z) = T _\mu \pi_C ( \nu)$, where $\mu =\JJ(\zt)\in
\mathfrak{g}^\ast  $ and $\nu = T_{\zt}\JJ(v_z)\in \mathfrak{g}^\ast$.

Consequently, $T _z\mathbf{K}(v_z)= 0$ is equivalent to
$T_{\zt}\JJ(v_z)\in\mathrm{Lie}(\Hbar)\subset\Hbar$, or equivalently
${\bf i}_{v_z}\omega\in\mathrm{Lie}(\Hbar)$, so that
$$\ker T_z\KK = \left[\left(\mathrm{Lie}(\Hbar)\right)^\circ\cdot z\right]^\omega.$$
Here ${\rm Lie}(\overline{\mathcal{H}}) \subset \mathfrak{g}^\ast$
is the Lie algebra of $\overline{\mathcal{H}}$, and ${\rm
Lie}(\overline{\mathcal{H}})^\circ$ its annihilator in $\gg$, and
the upper index $\omega$ denotes the $\omega$-orthogonal complement
of the set in question. The notation $\mathfrak{k} \cdot m $ for any
subspace $\mathfrak{k} \subset \mathfrak{g}$ has the usual meaning:
namely the vector subspace of $T _zM $ formed by evaluating all
infinitesimal generators $\eta_M$ at the point $z \in M $ for all
$\eta \in \mathfrak{k}$. Furthermore, ${\rm range}\, (T _z
\mathbf{K})= T _\mu \pi_C \left((\mathfrak{g}_{z})^\circ \right)$
(the Bifurcation Lemma).

\paragraph{Equivariance properties of the cylinder valued
momentum map.}  There is a $G$-action on $\gg^*/\Hbar$ with respect
to which the cylinder valued momentum map is $G$-equivariant. This
action is constructed by noticing first that since $G$ is connected
it follows (see~\cite{OR06}) that the Hamiltonian holonomy ${\cal
H}$ is pointwise fixed by the coadjoint action, that is, $ \mbox{\rm
Ad}^\ast _{g ^{-1}} h=h$, for any $g \in G $ and any $h\in\mathcal
H$.  Hence, the coadjoint action on $\gg^*$ descends to a well
defined action $\mathcal{A}d^\ast$ on $\gg^*/\Hbar$ defined so that
for any $g \in G $, $\mathcal{A}d^\ast_{ g ^{-1}} \circ \pi _C= \pi
_C \circ \mbox{\rm Ad} ^\ast  _{g ^{-1}}$. With this in mind, we
define
$\overline{\sigma}_{\mathbf{K}}: G\times
M \rightarrow \mathfrak{g}^\ast /\overline{\mathcal{H}}$ by
\[
\label{definition of the non-equivariance cocycle}
\overline{\sigma}_{\mathbf{K}}(g,z):= \mathbf{K}(g\cdot z) -
\mathcal{A}d^\ast_{g ^{-1}} \mathbf{K} (z).
\]
Since $M $ is connected by hypothesis, it can be shown that  $
\overline{\sigma}_{\mathbf{K}} $ does not depend on the point $z \in  M$ and
hence it defines a map $\sigma_{\mathbf{K}}: G \rightarrow  \mathfrak{g}^\ast
/\overline{\mathcal{H}} $ which is a group valued one-cocycle: for
any $g,h \in G $, it satisfies the equality $ \sigma_{\mathbf{K}}(gh)=
\sigma_{\mathbf{K}}(g)+
\mathcal{A}d^\ast _{g ^{-1}} \sigma_{\mathbf{K}} (h)$. This guarantees that the
map
\[
\begin{array}{cccc}
\Phi:& G \times  \mathfrak{g}^\ast
/\overline{\mathcal{H}}&\longrightarrow &
\mathfrak{g}^\ast /\overline{\mathcal{H}}\\
   &(g,\, \pi_C(\mu))&\longmapsto&
\mathcal{A}d^\ast _{g ^{-1}}(\pi_C(\mu))+ \sigma_{\mathbf{K}} (g),
\end{array}
\]
defines a $G$-action on $\mathfrak{g}^\ast /\overline{{\mathcal H}}$
with respect to which the cylinder valued  momentum map $\mathbf{K}$
is $G$-equivariant; that is, for any $g  \in G $, $z \in M $, we
have
\[
\mathbf{K}(g\cdot z) = \Phi(g,\mathbf{K}(z)).
\]
We will refer to $\sigma_{\mathbf{K}}: G \rightarrow  \mathfrak{g}^\ast
/\overline{\mathcal{H}} $ as the \emph{non-equivariance one-cocycle}
of the cylinder valued momentum map $\mathbf{K}:M \rightarrow
\mathfrak{g}^\ast/ \overline{\mathcal{H}} $ and to $\Phi$ as the
\emph{affine $G$-action} on $\mathfrak{g}^\ast /\overline{{\mathcal
H}}$ induced by $\sigma_{\mathbf{K}}$. The infinitesimal generators of the affine
$G$-action on $\mathfrak{g}^\ast/\overline{\mathcal{H}}$ are given
by the expression
\begin{equation}
\label{infinitesimal generators of the affine}
\xi_{\mathfrak{g}^\ast/ \overline{\mathcal{H}}} (\pi _C (\mu))=-T
_\mu \pi _C \left(\Psi (z)(\xi, \cdot ) \right),
\end{equation}
for any $\xi \in  \mathfrak{g}$, where  $\KK(z)=\pi_C(\mu)$, and
$\Psi:M \rightarrow Z ^2(\mathfrak{g})$ is the Chu map defined in
(\ref{eq:chu}).

The non-equivariance cocycles $\sigma_{\mathbf{J}}:
\widetilde{G}\rightarrow \mathfrak{g}^\ast$ and
$\sigma_{\mathbf{K}}: G \rightarrow
\mathfrak{g}^\ast/\overline{{\mathcal H}} $ are related by
\begin{equation}
\label{relation cocycles up and down} \pi _C \circ
\sigma_{\mathbf{J}}= \sigma_{\mathbf{K}}\circ q_G.
\end{equation}

\begin{proposition}
If the action of $G$ has an isotropic orbit then the cylinder valued
momentum map for this action can be chosen coadjoint equivariant.
\end{proposition}

\begin{proof}
This follows from Remark \ref{rmk:isotropic}.  Let $z _0\in M $ be a
point in the isotropic orbit and construct a universal cover
$\widetilde{M}$ of $M$ by taking homotopies of curves with a fixed
endpoint starting at $z_0$.  Let $\mathbf{J}: \widetilde{M}
\rightarrow  \mathfrak{g}^\ast $ be the momentum map for the
$\widetilde{G}$-action on $\widetilde{M} $  introduced in
Proposition~\ref{momentum map for the lifted g tilde action}. Since
the $G$-orbit containing $z _0 $ is isotropic, the integrand
in~(\ref{non-equivariance cocycle}) is  identically zero and hence
$\sigma_{\mathbf{J}} =0$ (see Remark \ref{rmk:isotropic}). Therefore
by (\ref{relation cocycles up and down}) the non-equivariance
cocycle $\sigma_{\mathbf{K}}$ satisfies $\sigma_{\mathbf{K}}\circ
q_G = \pi_C \circ \sigma_{\mathbf{J}}=0 $.
\end{proof}

\begin{remark}
\label{smart choice of j} \normalfont For any Hamiltonian cover
$p _N:N \rightarrow  M $ of $(M, \omega)$ there exists a momentum
map $\mathbf{J}_N:N \rightarrow \mathfrak{g}^\ast$ for the
$\widetilde{G}$ (and also $G _N$) action on $N$ such that
$\mathbf{J}_N \circ q _N= \mathbf{J}$ and $\sigma_{\mathbf{J}_N}=
\sigma_{\mathbf{J}}$, where $q _N  : \widetilde{M}\rightarrow N$ is
the $\widetilde{G}$-equivariant cover such that  $p _N\circ
q_N=q_M$.  Consequently, there is a commutative diagram analogous to (\ref{diagram commutative cylinder valued momentum map}) with $N$ and $\JJ_N$ in place of $\Mt$ and $JJ$.
\end{remark}

%%%%%%%%%%%%%%%%%%%%%%%%%%
\subsection{Reductions}

The following result establishes a crucial relationship between the
deck transformation group of $q_M: \Mt\rightarrow M $, that is,
$\Gamma:=\pi_1(M, z _0)$, and the deck transformation group of
$\widehat{p}: \widehat{M} \rightarrow M $, that is ${\mathcal
H}\simeq \Gamma/\Gamma_0$.

\begin{proposition} \label{prop:relation holonomy proposition} Let $G$
be a connected Lie group acting symplectically on the symplectic
manifold $(M, \omega)$ with Hamiltonian holonomy ${\mathcal H}$ and
let $\JJ: \Mt\rightarrow M $ be the momentum map for the lifted
action on $(\Mt,\zt_0)$ defined in Proposition~\ref{momentum map for
the lifted g tilde action}. Then, for any $\mu \in
\mathfrak{g}^\ast$
\begin{equation} \label{eq:relation holonomy deck}
q_M^{-1}\left(q_M(\JJ^{-1}(\mu))\right) =  \JJ^{-1}(\mu+ \H).
\end{equation}
More generally, for any Hamiltonian cover $p_N:(N,y_0)
\rightarrow (M,z_0)$ of $(M, z_0,\omega)$, let $\JJ_N : N
\rightarrow \mathfrak{g}^\ast$ be the momentum map discussed in
Remark~\ref{smart choice of j}. Then, for any $\mu \in
\mathfrak{g}^\ast$
\begin{equation} \label{eq:relation holonomy deck general}
p_N^{-1}\left(p_N(\mathbf{J}_N^{-1}(\mu))\right) =
\mathbf{J}_N^{-1}(\mu + \mathcal{H}).
\end{equation}
\end{proposition}

\begin{proof}
Since $\Gamma$ acts transitively on the fibres of $q_M$,
(\ref{eq:relation holonomy deck}) is equivalent to
$$\JJ^{-1}(\mu+ \H) = \Gamma\cdot \JJ^{-1}(\mu).$$
By Proposition \ref{momentum map for the lifted g tilde action}, if
$\JJ(\zt)=\mu$ and $\gamma\in\Gamma$ then $\JJ(\gamma\cdot\zt) =
\mu+\nu$ for some $\nu\in\H$; that is,
$\gamma\cdot\zt\in\JJ^{-1}(\mu+\H)$.  Conversely, given $\nu\in\H$
there is a $\gamma\in\Gamma$ for which $\JJ(\gamma\cdot\zt) =
\mu+\nu$ so proving the statement.

In order to prove~(\ref{eq:relation holonomy deck general}) let $q
_N : \widetilde{M}\rightarrow N$ be the $\widetilde{G}$-equivariant
cover such that  $p _N\circ q _N = q_M$. This equality and the
surjectivity of $q _N $ imply that for any set $A \subset N $, $p _N
(A)= q_M (q _N ^{-1} (A)) $. Now, the relations $\mathbf{J}_N \circ
q _N = \mathbf{J}  $ and~(\ref{eq:relation holonomy deck}) imply
that $ q_M \left( q _N ^{-1} (\mathbf{J} _N ^{-1}(\mu + {\mathcal
H}))\right)= q_M\left(q _N ^{-1}(\mathbf{J}_N ^{-1} (\mu)) \right) $
and hence $ p_{N}(\mathbf{J}_N^{-1}(\mu+ {\mathcal
H}))=p_{N}(\mathbf{J}_N^{-1}(\mu))$, as required.
\end{proof}

The main result of this section shows that when the Hamiltonian holonomy is closed
reduction  behaves well with respect to the lifting of the action to
any Hamiltonian cover. More explicitly, we show that in order to
carry out reduction one can either stay in the original manifold and
use the cylinder valued momentum map or one can lift the action to a
Hamiltonian cover, perform ordinary symplectic (Marsden-Weinstein)
reduction there and then project the resulting quotient. The two
strategies yield closely related results.  Notice that if the
Hamiltonian holonomy of the action ${\mathcal H} $ is not closed in
$\mathfrak{g}^\ast$, the reduced spaces obtained via the cylinder
valued momentum map are in general not symplectic but Poisson
manifolds~\cite{OR06}.

For the remainder of this section we assume the Hamiltonian holonomy
$\H$ to be a closed subset of $\gg^*$, and we write $\gt\cdot\mu$
for the \emph{modified} coadjoint action of $G'$ or $\Gt$ on $\gg^*$,
and similarly $g\cdot[\mu]$ for the inherited action on $\gg^*/\H$.
We also write $\Gamma':=\image(a_{z_0})$, where $a_{z_0}$ is defined in
(\ref{eq:a_z_0}).

Let $N$ be any Hamiltonian cover of $M$, and consider the diagram for $N$ analogous to (\ref{diagram commutative cylinder valued momentum map}); of
course particular cases of interest are $N=\Mt$ and $N=\Mhat$. As
$\H$ is closed, the image of $\JJ_N^{-1}(\mu+\H)$ under $p_N$ is
precisely $\KK^{-1}([\mu])$, by the definition of $\KK$. Reduction
of each defines a map
$$(p_N)_\mu:N_\mu \longto M_{[\mu]}.$$
In the case that $N=\Mt$, we denote the projection by
$(q_M)_\mu:\Mt_\mu\to M_{[\mu]}$.

For each $\mu\in\gg^*$ define
$$\Gamma_\mu = \Gamma \cap \JJ^{-1}(\sigma_\mu(\Gt))$$
where $\sigma_\mu:\Gt\to\gg^*$ is the 1-cocycle $\sigma_\mu =
\sigma_\JJ + \delta\mu$ and $\delta\mu(\gt) = \delta\mu(g)=
\Ad^*_{g^{-1}}\mu-\mu$ is the coboundary associated to $\mu$.
Note that for all $\mu\in\gg^*$, $\Gamma' < \Gamma_\mu$.
Indeed, given $\gt\in\pi_1(G,e)$, $\JJ(\gt\cdot\zt_0) = \sigma(\gt) = \sigma_\mu(\gt)$ as required; the last equality holds because for
$\gt\in\pi_1(G,e)$, $\delta\mu(\gt) = 0$.

Furthermore, we have that $\Gamma_\mu\supset\Gamma_0 = \JJ^{-1}(0)\cap\Gamma$. Since both $\Gamma\,'$ and $\Gamma_0$ are normal subgroups of $\Gamma$
(and hence of $\Gamma_\mu$), with $\Gamma\,'$ being in the
centre, it follows that, for all $\mu\in\gg^*$, the product
\begin{equation}\label{eq: Gamma_mu}
\Gamma\,'\Gamma_0 \lhd \Gamma_\mu.
\end{equation}

\begin{theorem} \label{thm:reduced covering} Suppose the action
of $G$ on $(M,\omega)$ is free and proper, and the holonomy group
$\H$ is closed.  Then the map $(q_M)_\mu:\Mt_\mu\to M_{[\mu]}$ is a
cover, with transitive deck transformation group isomorphic to
$$\Gamma_{\mu,\mathrm{red}} := \Gamma_\mu/\Gamma\,'.$$
More generally, if $N$ is a normal Hamiltonian cover of $M$ then
$(p_N)_\mu$ is a normal cover, with the deck transformation group
$$\Gamma_\mu/\left(\Gamma_N\Gamma\,'\right).$$
\end{theorem}

\begin{proof}
We approach this from the point of view of orbit reduction; that is
we consider
$$M_{[\mu]} = \KK^{-1}(G\cdot[\mu])/G\subset M/G,\quad
\mbox{and}\quad \Mt_\mu = \JJ^{-1}(\Gt\cdot\mu)/\Gt \subset \Mt/\Gt.
$$
In both cases, the $G$ or $\Gt$ actions are the coadjoint action
modified by the cocycle $\sigma_\KK$ and $\sigma_\JJ$, respectively.
It is well-known that for proper actions, point and orbit reductions
are equivalent (for a proof, see Theorem 6.4.1 of \cite{hsr}), and
the equivalence respects the projections induced by $\Mt\to M$.

Consider then the following commutative diagrams:

\begin{equation} \label{eq:reduction diagram}
\begin{CD}
\JJ^{-1}(\Gt\cdot\mu)@>\pi_{\Mt}>> \Mt_\mu \\
@V{q_M}VV @VV{q_M'}V\\
\KK^{-1}(G\cdot[\mu])@>\pi_M>>M_{[\mu]}\\[6pt]
\end{CD}
\qquad\quad\subset\qquad\quad
\begin{CD}
\Mt@>\pi_{\Mt}>>  \Mt/G' \\
@V{q_M}VV @VV{q_M'}V\\
 M @>\pi_M>> M/G\\[6pt]
\end{CD}
\end{equation}
The maps in the left-hand diagram are just restrictions of those in
the right-hand one.

First we claim that $q_M: \JJ^{-1}(\Gt\cdot\mu) \to \KK^{-1}(G\cdot[\mu])$ is a cover whose group of covering transformations is $\Gamma_\mu$ defined above. The result then follows from Proposition \ref{prop:orbit space covering}, but with $\Gamma$ replaced by $\Gamma_\mu$, since $\Gamma\,' < \Gamma_\mu$.

To prove the claim, we know from Proposition \ref{prop:relation
holonomy proposition} that $q_M^{-1}(\KK^{-1}([\mu])) =
\JJ^{-1}(\mu+\H)$. Saturating by $\Gt$, we have
$$q_M^{-1}(\KK^{-1}(G\cdot[\mu])) = \JJ^{-1}(\Gt\cdot(\mu+\H)),
$$
and this is a cover with group $\Gamma$ (that of the cover
$\Mt\to M$).

Now let $z\in M$ be such that $\KK(z)=[\mu]$ (so in particular
$z\in\KK^{-1}(G\cdot[\mu])$), and let $Z=q_m^{-1}(z)$ be the fibre
over $z$. If $\zt\in Z$ then $Z=\Gamma\cdot\zt$, and
$\JJ(\Gamma\cdot \zt) = \mu+\H$, so we choose $\zt\in Z$ such that
$\JJ(\zt)=\mu$.

We now show that $Z\cap\JJ^{-1}(\Gt\cdot\mu)=\Gamma_\mu\cdot\zt$. To
this end, let $\zt_1\in Z$. Then $\exists \gamma\in\Gamma$ such that
$\zt_1=\gamma\cdot\zt$, so
$$\JJ(\zt_1)=\JJ(\zt)+\JJ(\gamma) = \mu+\JJ(\gamma).
$$
Then $\mu+\JJ(\gamma)\in \Gt\cdot\mu$ if and only if
$\exists\,\gt\in\Gt$ such that
$$\mu+\JJ(\gamma) = \gt\cdot\mu = \Ad^*_{g^{-1}}\mu+\sigma(\gt),$$
so that $\JJ(\gamma)= \delta\mu(\gt) + \sigma(\gt) =
\sigma_\mu(\gt)$; that is, $\gamma\in\Gamma_\mu$, as required.

The proof of the second part of the theorem, with a general normal
cover $N$, is identical, given that $N=\Mt/\Gamma_N$.
\end{proof}

\begin{corollary}
The cover $\Mhat_\mu \to M_{[\mu]}$ has cover transformation
group $\Gamma_\mu/\Gamma_0\Gamma\,'$. This is trivial if
$\JJ(\Gamma\,') = \H\cap \sigma_\JJ(\Gt)$, in which case the
cover is a symplectomorphism.
\end{corollary}

\begin{remark}
If the Hamiltonian holonomy is not closed but the action is still
free and proper, the reduced spaces $M_{[\mu]}$ and $\Mt_\mu$ are
Poisson manifolds \cite{OR06}, and the natural map $p_\mu:\Mt_\mu\to M_{[\mu]}$
is a surjective Poisson submersion.
\end{remark}

%%%%%%%%%%%%%%%%%%
\subsection{Example}
We continue the example of $G$ acting on $T^*G$ with symplectic form modified by a cocycle $\theta$, as discussed in \S\ref{sec:example of T^*G}. In this case, $\Gamma = \pi_1(G,e)$ and $a_{z_0}:\pi_1(G,e) \to \Gamma$ is the identity, so $\Gamma\,'=\Gamma$ and it follows that $\Gamma_\mu=\Gamma$ for all $\mu\in\gg^*$.

Write $M=T^*G$ and $\Mt=T^*\Gt$ and assume that the Hamiltonian holonomy $\H=\Theta(\Gamma)\subset\gg^*$ is closed. It follows from Theorem \ref{thm:reduced covering} that the projection $\Mt_\mu \to M_{[\mu]}$ is a cover with trivial (and transitive) deck transformation group, so is in fact a symplectomorphism.  Indeed the same is true for any intermediate cover $G_1$ for which the action on $T^*G_1$ is Hamiltonian.
In particular, we find that for the left action of $G$ on $T^*G$ with modified symplectic form, Hamiltonian reduction for a Hamiltonian lift and symplectic reduction via the cylinder valued momentum map yield the same result.

The well-known statement that the symplectic reduced spaces for the canonical left action of $G$ on $T^*G$ coincide with the coadjoint orbits \cite{mwr} remains true when both the symplectic structure and the action on $\gg^*$ are modified by a cocycle $\Theta$ (see for example \cite{hsr}).  The statement above shows that this remains true for cylinder valued momentum maps, where the orbits are those of $\Gt$ in $\gg^*$ rather than those of $G$ in $C$.

\begin{example} Returning to Example \ref{eg:torus} on the torus, given $\theta\in H^1_s(\tt,\tt^*)$ the orbits of the modified coadjoint action of $\RR^d$ are the affine subspaces parallel to $\image(\theta)\subset\tt^*$, and so the reduced spaces for this action are symplectomorphic to these affine subspaces.  If $\theta$ is chosen so that the holonomy is closed (eg, $d$ is even and $\theta$ is invertible) then the same is true of the reduced spaces for the action of $\TT^d$ on $T^*\TT^d$ via the cylinder valued momentum map.

\end{example}

\begin{example}
Returning now to Example \ref{eg:Heisenberg}, the symplectic reduced spaces for the Heisenberg group with the symplectic structure $\Omega_{\mathrm{canon}}+\pi^*B_\Sigma$ on $T^*H$ are the orbits for the modified coadjoint action.  Calculations show these to be the level sets of the Casimir function $f(\psi,\nu) = \half\psi^2 - \omega^{-1}(\sigma,\nu)$, which are parabolic cylinders.  Since the Hamiltonian holonomy $\H$ is closed, it follows from the results above that the same is true for reduction via the cylinder valued momentum map on $T^*G$.
\end{example}

\paragraph{Acknowledgments.} The authors would like to thank Liviu Ornea for enlightening discussions that motivated part of this work, and Rui Loja Fernandes for explaining the relationship with the groupoids approach of Mikami and Weinstein. JPO thanks the University of Manchester for its hospitality while part of this work was carried out. JPO has been partially supported by the French Agence National de la Recherche,
contract number JC05-41465 and by a ``Bonus Qualit\'e Recherche"
contract from the Universit\'e de Franche-Comt\'e.

%%%%%%%%%%%%%%%%%%%%%%%%%%%%%%%%%%
%%%%%%%%%%%%%%%%%%%%%%%%%%%%%%%%%%

\small

\bigskip

\hbox to 0.8\textwidth{\hrulefill}

\smallskip

\parindent=0pt

\parbox[t]{0.4\textwidth}{{\sl School of Mathematics, \\
University of Manchester, \\
Oxford Road,\\
Manchester M13 9PL, \\
UK.}\\

\texttt{j.montaldi@manchester.ac.uk} }
\parbox[t]{0.5\textwidth}{{\sl Centre National de la Recherche Scientifique, \\
D\'epartement de Math\'ematiques de Besan\c{c}on,\\
Universit\'e de Franche-Comt\'e, \\
UFR des Sciences et Techniques, \\
16 route de Gray, \\
25030 Besan\c{c}on c\'edex,\\
France.}\\

\texttt{Juan-Pablo.Ortega@univ-fcomte.fr} }

%%%%%%%%%%%%%%%%%%%%%%%%%%%%%%%%%%%%%%
\end{document}